\documentclass[12pt,twoside,reqno]{amsart}

\usepackage{amsmath,amssymb,amsthm}
\usepackage{fancyhdr}
\usepackage{epsfig,color,colordvi} 
\usepackage[numbers,sort&compress,square,comma]{natbib}

\theoremstyle{definition}

\newtheorem{example}{\sc Example}[section]

\theoremstyle{remark}
\newtheorem{remark}{\sc Remark}[section]

\theoremstyle{plain}
\newtheorem{theorem}{\sc Theorem}[section]
\newtheorem{lemma}{\sc Lemma}[section]
\newtheorem{proposition}{\sc Proposition}[section]

\newcommand*{\be}{\begin{equation}}
\newcommand*{\ee}{\end{equation}}
\newcommand*{\ba}{\begin{aligned}}
\newcommand*{\ea}{\end{aligned}}
\newcommand*{\nn}{\nonumber}

\def\V{\bbV}
\def\vtil{{\tilde{v}}}

\def\kD{\mathfrak{D}}
\def\Qtil{{\widetilde Q}}
\def\pbar{{\bar{p}}}
\def\pibar{{\bar{\pi}}}
\def\zbar{{\bar{z}}}

\def\Btil{{\widetilde{B}}}

\def\Stil{{\widetilde{S}}}
\def\Ytil{{\widetilde{Y}}}

\def\e{\varepsilon}
\def\ind{\mathbf{1}}  % Careful that Firas also has \one 
\def\cJ{\mathcal{J}}  % and Pb (same as bbP!)
\def\cG{\mathcal{G}}

\newcommand*{\si}{\sigma}

\newcommand*{\la}{\lambda}

\def\cB{\mathcal{B}}
\def\cC{\mathcal{C}}
\def\uhat{{\hat u}}

\def\pihat{{\hat\pi}}

\def\kS{{\mathfrak S}}
\def\Xhat{{\widehat{X}}}

\def\bbR{\mathbb R}
\def\bbN{\mathbb N}
\def\bbZ{\mathbb Z}
\def\bbQ{\mathbb Q}
\def\bbP{\mathbb P}
\def\bbE{\mathbb E}
\def\bbV{\mathbb V}

\def\fQ{{{\rm Q}\kern-.65em {}^{{}_/ }\,}}
\def\fQQ{ {{\rm Q}\kern-.57em \scriptscriptstyle{}^{]\kern.055em[}\,}}

\def\one{{{\rm 1\mkern-1.5mu}\!{\rm I}}}
\def\w{\omega}

\DeclareMathOperator{\Var}{Var}
\def\VVar{\text{$\mathbb{V}$\!ar}}

\def\lsup{\varlimsup}

\def\ord{\kern0.1em o\kern-0.02em{}_{\ds\breve{}}\kern0.1em}
\def\Ord{\kern0.1em O\kern-0.02em{\ds\breve{}}\kern0.1em}
\def\ds{\displaystyle}

\def\fmonth{\ifcase\month\or Jan\or Feb\or Mar\or Apr
\or May\or Jun\or Jul\or Aug\or Sep
\or Oct\or Nov\or Dec\fi\ }
\def\mmddyyyy{\the\month.\the\day.\the\year}
\def\ddmmyyyy{\the\day.\the\month.\the\year}
\def\Mddyyyy{\fmonth~\the\day,~\the\year}

\def\R{\bbR}
\def\N{\bbN}
\def\Z{\bbZ}
\def\Q{\bbQ}
\def\P{\bbP}
\def\E{\bbE}

\providecommand{\abs}[1]{\left\vert#1\right\vert}
\providecommand{\norm}[1]{\left\Vert#1\right\Vert}

\numberwithin{equation}{section}

\allowdisplaybreaks[1]

\begin{document}

%% Please copy the part bellow into the textfile

%%%%%%%
%%% Author(s) and other info
%%%%%%%

\author{Firas Rassoul-Agha}
\address{F.\ Rassoul-Agha, Department of Mathematics, University of Utah,
155 South 1400 East, Salt Lake City, UT 84112-0090, USA.}
%\curraddr{}
\email{firas@math.utah.edu}
%\urladdr{www.smith.com}
%\thanks{NSF}
\author{Timo Sepp\"al\"ainen}
\address{T.\ Sepp\"al\"ainen, Mathematics Department, University of Wisconsin-Madison,
Van Vleck Hall, 480 Lincoln Dr, Madison WI 53706-1388, USA.}
%\curraddr{}
\email{seppalai@math.wisc.edu}
%\urladdr{www.doe.co.uk}
\thanks{T.~Sepp\"al\"ainen was partially supported by
National Science Foundation grant DMS-0402231.}

\date{\today}
% this will put the date as a footnote on the first page

%\translator{Nobody}
%\dedicatory{Joe}

%%%%%%
%%% If you want an abstract
%%%%%%

%% Uncomment this if using \bfsection
%\renewcommand\abstractname{\noindent\bf Abstract}

\frenchspacing

\begin{abstract}
We consider a ballistic random walk in an i.i.d. random environment
that does not allow retreating in a certain fixed direction.
 Homogenization and regeneration techniques  combine to
prove a law of large numbers and an averaged invariance principle.
The assumptions are non-nestling and
$1+\e$ (resp.\ $2+\e$) moments for the step  of the
walk uniformly in the environment, for the law of large numbers
(resp.\ invariance principles). We also investigate 
 invariance
principles under fixed environments, and invariance
principles for the environment-dependent  mean of the walk. 
\end{abstract}

%%%%%%%
%%% Don't forget the title and date. To keep track of the notes.
%%% Preferably, name the file with something indicating the
%%% title and the date, as well.
%%%%%%%

%% You can copy the part bellow into mytext

%% STARTING HERE

\title[Random Walk in Random Environment]
{Ballistic Random Walk in a Random Environment with a Forbidden
Direction}

\maketitle

\def\R{\bbR}
\def\N{\bbN}
\def\Z{\bbZ}
\def\Q{\bbQ}
\def\P{\bbP}
\def\E{\bbE}
\def\V{\bbV}

%\frenchspacing

\section{Introduction}
\label{intro}
This paper studies random walk in a random
environment (RWRE)  on the  $d$-dimensional integer lattice $\Z^d$. This
is a basic model in the field of disordered or random media. The
dimension $d$ is in general
 any positive integer, but in certain results we are forced
to  distinguish between $d=1$ and $d\geq 2$.

A general description of  this model follows.
An environment is a configuration of  vectors of jump probabilities
\[\w=(\w_x)_{x\in\Z^d}\in\Omega={\mathcal P}^{\Z^d},\] where
$\mathcal P=\{(p_z)_{z\in\Z^d}:p_z\geq0,\sum_z p_z=1\}$ is the simplex of all
(infinite) probability vectors on $\Z^d$.
We use the notation $\w_x=(\pi_{x,x+y})_{y\in\Z^d}$ for the coordinates of $\w_x$.
%$\Omega$  denotes the space of all such transition probabilities.
The space $\Omega$ is equipped with the canonical product
$\si$-field $\kS$ and with the natural shift
$\pi_{xy}(T_z\w)=\pi_{x+z,y+z}(\w)$, for $z\in\Z^d$.
On the space of environments $(\Omega,\kS)$ we are given a
$T$-invariant probability measure $\P$ such that the system
$(\Omega,\kS,(T_z)_{z\in\Z^d},\P)$ is ergodic.
%We will denote expectation under $\P$ by $\E$.
All our results are for  i.i.d.\  environments.
This means that the random probability
 vectors $(\w_x)_x$ are i.i.d.\ across the sites
$x$ under $\P$,
or in other words  that $\P$ is an i.i.d.\  product measure on ${\mathcal P}^{\Z^d}$.

Here is how the random walk operates.
First an environment $\w$ is chosen
from the distribution $\P$. The environment $\w$ remains fixed for
all time. Pick an initial state $z\in\Z^d$.
 The random walk in environment $\w$ started at $z$ is then the canonical
Markov chain $\Xhat=(X_n)_{n\geq0}$ with state space $\Z^d$
whose path measure  $P_z^\w$ satisfies
 \begin{align*}
P_z^\w(X_0=z)&=1 & &\text{(initial state)},\\
P_z^\w(X_{n+1}=y|X_n=x)&=\pi_{xy}(\w)
& &\text{(transition probability).}
\end{align*}
The probability distribution  $P_z^\w$ on random walk paths
 is called the {\sl quenched law}. %The corresponding expectation is $E_z^\w$.
The joint  probability distribution
 \[P_z(d\Xhat,d\w)=P_z^\w(d\Xhat)\P(d\w)\]
on walks and environments is called the
 {\sl joint annealed law}, while its marginal on walks
  $P_z(d\Xhat,\Omega)$ is called simply
  the {\sl annealed law}.
We will use $\E$, $E_z^\w$, and $E_z$ to denote expectations under, respectively,
$\P$, $P_z^\w$, and $P_z$.

Not much of interest has been proved for random walk in random
environment in this generality, except in one dimension.
In this paper and the subsequent paper \cite{forbidden-qclt} 
we  study the invariant distributions, 
laws of large numbers and central limit theorems for a 
class of ballistic walks in an i.i.d.~environment 
 that possess a ``forbidden direction.''  
More precisely, we impose a non-nestling assumption that creates a
drift in some spatial direction $\uhat$, and we also
 prohibit the walk from retreating in
direction $\uhat$. However, there is some freedom in the
choice of $\uhat$. The long term velocity $v$ of the walk
need not be in direction $\uhat$, although of course the
assumptions will imply $\uhat\cdot v> 0$. These assumptions are 
introduced in Section \ref{assumptions}.  Given the hypotheses
described above our results are
complete, except for some moment assumptions that most
certainly are not best possible.  

Our main result is the central limit  picture for this class of
walks. As preparation for this we first need to establish the 
existence of suitable invariant distributions for the environment
chain, and also prove a law of large numbers $n^{-1}X_n\to v$. 

The first point of the fluctuation 
question is that  the annealed invariance principle (functional
 central limit theorem) 
 with $nv$ centering
is always valid. In other words, the distribution
 of the process $\{B_n(t)=n^{-1/2}(X_{[nt]}-[nt]v):  t\geq0\}$
under $P_0$ converges to a Brownian motion with a certain diffusion matrix
that we give a formula for. This will be the object of Section
\ref{inv1}.

The second fluctuation issue is the quenched invariance principle
under measures $P^\w_0$. 
Here we encounter the following dichotomy. 

\smallskip

I. Suppose the walk is one-dimensional, or that the environment $\w$ 
constrains the walk to be essentially one-dimensional in a sense
made precise by  Hypothesis (E) on page \pageref{hypoE} below. Then 
process $B_n$ does not satisfy a quenched invariance principle.
However, under quenched mean centering a quenched invariance principle 
does hold. In other words,  for $\P$-almost every $\w$, 
 the process $\{\Btil_n(t)=n^{-1/2}(X_{[nt]}-E^\w_0(X_{[nt]})):  t\geq0\}$
under $P^\w_0$ 
converges weakly to a certain Brownian motion. Furthermore, 
the centered quenched mean process 
$\{n^{-1/2}(E^\w_0(X_{[nt]})-ntv):  t\geq0\}$ also satisfies 
an invariance principle. 

\smallskip

II. In the complementary case the walk explores its environment
thoroughly enough to suppress the fluctuations of the quenched mean. 
The variance \[\E[\,\lvert E^\w_0(X_n)-E_0(X_n)\rvert^2\,]\]
of the quenched mean is of order $n^{1-\delta}$ for $\delta>0$, the centerings 
$nv$ and $E^\w_0(X_n)$ for the walk are asymptotically indistinguishable
on the diffusive scale, and
both processes $B_n$ and $\Btil_n$ defined above satisfy  a quenched
invariance principle with the same limiting Brownian motion. 

\smallskip

We can also view the dichotomy through the 
following  decomposition of the fluctuations:
%of the annealed covariance of the walk 
% $X_n$ into the covariance of the quenched
%mean plus the mean of the quenched covariance:
\begin{align}
&E_0[(X_n-nv)(X_n-nv)^t]\nonumber\\
&\quad=
\E[E_0^\w(X_n X_n^t)-E_0^\w(X_n)E_0^\w(X_n)^t]+\E[(E_0^\w(X_n)-nv)(E_0^\w(X_n)-nv)^t]. 
\label{decomposition}
\end{align}
In case I both terms contribute on the diffusive scale all the way
to the limit. As the results in  Section \ref{inv2} show, 
the fluctuations described by 
the annealed central limit theorem decompose  into quenched fluctuations
of the walk around  the quenched mean, plus 
the fluctuations of  the quenched mean around $nv$. (The
 centering $nv$ in  \eqref{decomposition}  is not exactly the true mean $E_0(X_n)$, but
the difference of the two is bounded as shown in Theorem \ref{lln}.)

In case II the second term 
on the right hand side of  \eqref{decomposition} is
 negligible relative  to the two other terms of order $n$. On the 
diffusive scale the fluctuations of the quenched mean are not felt. 

This paper covers the type I fluctuation results in Section \ref{inv2}
which forms the main portion of the paper. The type II
results are reported in the companion paper  \cite{forbidden-qclt}.
An appendix in the present paper contains a renewal lemma which
is used in both papers. 
Despite our restriction to the forbidden direction
case,  we cautiously believe that
the results and methods of proof are characteristic of more general ballistic
random walks in random environment.  
%In addition to the quenched 
%invariance principle itself,  an interesting question for future 
%work is whether the variance of the quenched mean is of order $n^{1/2}$
%in more general ballistic walks. 

%In all other cases the quenched mean has variance of order $n^{1/2}$,
%so in the limit this term does not contribute to the variance
%decomposition.  Thus on the diffusive scale the quenched walk
%will obey the same central limit theorem regardless of the centering chosen,
%$nv$ or $E^\w_0(X_n)$. 

We turn to a brief overview of some past work. 
Two main approaches to laws of large numbers and central limit 
theorems  for %random walk in random environment 
RWRE stand out in the literature. 

Using regeneration techniques (see the renewal times $\si_k$ in \eqref{sik} below),
\cite{czlln} proves a law of large numbers,
\cite{czclt} proves an annealed central limit theorem, both in the case of a mixing
environment, while \cite{effective} proves both results when the environment is product.
Also, using these techniques, \cite{01llniid}
reduces the proof of the law of large numbers in product environments to that of a
certain 0-1 law.

On the other hand, using homogenization techniques (see Section \ref{invariant}),
\cite{llnrwre} proves a law of large numbers in mixing environments and
\cite{qclt-spacetime} proves a quenched central limit theorem in space-time
product environments.
Furthermore, using this approach, \cite{01lln} extends the result of \cite{01llniid}
to mixing environments.

The majority of results have been shown for
the ballistic case where  the velocity of escape
does  not vanish. Recent quenched central limit results
 in the diffusive case can be found
in \cite{cutpoints} and \cite{cltrwrenoise}.
For general overviews of the field and further literature
we refer the reader to the lectures
 in \cite{stflour} and
\cite{rwrelectures}.
The introduction of \cite{qclt-spacetime} also presents
a brief list  of papers on
central limit theorems for random walks in random environment.

To prove our central limit results 
 we combine both homogenization and regeneration
techniques.  As groundwork  we need suitable
 equilibrium measures and the law of large numbers.
These have been
proved in the past for more general ballistic walks, but typically with
 stringent assumptions of ellipticity and nearest-neighbor steps.
Therefore, in Section \ref{invariant},
we reprove these results under our weaker assumptions.  
In particular, these preliminary results are then used in the 
paper  \cite{forbidden-qclt} that proves the quenched invariance
principle for the genuinely multidimensional walks.

%Here is how this article is organized. First, in Section \ref{assumptions},
%we introduce our assumptions. Then, in Section \ref{invariant}, we prove
%the law of large numbers. In Section \ref{inv1}, we prove the annealed
%invariance principle. To keep the article at a reasonable size, the quenched
%invariance principle under the ellipticity hypothesis (E) (see page \pageref{hypoE}
%below) is proved in a separate paper \cite{forbidden-qclt}.
%Finally, in Section \ref{inv2}, we show that the this ellipticity hypothesis
%is also necessary for the
%annealed invariance principle to generalize to a quenched one.

\section{Assumptions}
\label{assumptions}
The type of random walk in random environment that we
study in this article is defined by  the following asymmetry condition. Given a nonzero
vector $\uhat$, say the distribution
$\P$  on environments has {\sl forbidden direction $-\uhat$} if
\begin{align}
\P\left(\,
\sum_{z:z\cdot\uhat\geq 0}\pi_{0z}(\w)=1
\right) =1.
\label{forbidden}
\end{align}
This  condition says that  $X_n\cdot\uhat$ never
decreases along the walk.

Assumption \eqref{forbidden} alone does not give the regeneration
structure we need. Trivially, even $\P(\sum_{z\cdot\uhat>0}\pi_{0z}=0)=1$
is possible, which permits an arbitrary
$(d-1)$-dimensional walk and the forbidden direction condition is vacuous.
However, if we just assume $\P(\sum_{z\cdot\uhat>0}\pi_{0z}>0)>0$, then
there could still be a positive chance of having $x\not=y\in\Z^d$ such that
$\pi_{0x}+\pi_{0y}=\pi_{x0}=\pi_{y0}=1$. In this case, a walker starting at 0
will be stuck in the set $\{0,x,y\}$ forever and
the hypothesis on the $\uhat$-direction is again
useless. If we add on the assumption  $\P(\sum_{z\cdot\uhat>0}\pi_{0z}>0)=1$, then
the walker will eventually move in the
$\uhat$-direction, but the expected time of
doing so could be infinite, leading to zero
 velocity in direction $\uhat$. This scenario is illustrated by the next example. 

\begin{example}
\label{si=infty}
Consider a two-dimensional product environment $\P$ with marginal at $0$
given by
\[\P(\w_0=p_i)=\P(\w_0=q_i)=2^{-i},\]
where for $i\geq 1$
\[p_i=(1-\alpha_i)\delta_{(1,0)} + \alpha_i \delta_{(0,1)}\text{ and }
q_i=(1-\alpha_i)\delta_{(1,0)} + \alpha_i \delta_{(0,-1)}\]
with $\alpha_i = 1 - 2^{-2i}$.  The formulas above mean that
if $\w_x=p_i$ then from $x$ the walk jumps to $x+(1,0)$ with
probability $1-\alpha_i$, and to  $x+(0,1)$ with
probability $\alpha_i$. Note that the jump size is finite, which is even
stronger than the moment hypothesis (M). Also, assumption \eqref{forbidden}
holds with $\uhat=(1,0)$.

Define $\si_1=\inf\{n:X_n\cdot\uhat\geq X_0\cdot\uhat+1\}$. Then
\[P_0(\sigma_1>n)\geq
\sum_{i\geq 1} \P(\w_{(0,0)}=p_i,\w_{(0,1)}=q_i)\alpha_i^n
%[explanation: walker jumps n times between sites (0,0) and (0,1)]
= \sum_{i\geq 1} 2^{-2i} (1-2^{-2i})^n\]
and
\[E_0(\sigma_1)=\sum_{n\geq 0} P_0(\sigma_1>n)
\geq \sum_{i\geq 1} 2^{-2i} \sum_{n\geq 0} (1-2^{-2i})^n
= \sum_{i\geq 1} 2^{-2i} 2^{2i} = \infty.\]
\end{example}

These degenerate situations are prevented by a 
 {\sl non-nestling} condition  in direction $\uhat$. That is
our second assumption.

\newtheorem*{hypothesisN}{\sc Hypothesis (N)}
\begin{hypothesisN}
There exists a positive deterministic constant $\delta$ such that
\[\P\left(\sum_z z\cdot\uhat\,\pi_{0z}\geq\delta\right)=1.\]
\end{hypothesisN}

It is noteworthy at this point that with a forbidden direction $-\uhat$,
any sensible uniform ellipticity assumption implies hypothesis (N).
 Indeed, if there exists a $\kappa>0$ and an $x_0\in\Z^d$ with
$x_0\cdot\uhat>0$ and $\P(\pi_{0x_0}\geq\kappa)=1$, then (N) is satisfied
with $\delta=\kappa\, x_0\cdot\uhat>0$.  In this sense, hypothesis (N) is
a fairly weak assumption that insures a non-zero velocity in direction
$\uhat$.

Most of the literature on random walks in random environment concentrates
on walks with nearest-neighbor  or bounded jumps.
This restriction we do not need, but moment assumptions are necessary
for laws of large numbers and central limit theorems. This is our next
assumption.

\newtheorem*{hypothesisM}{\sc Hypothesis (M)}
\begin{hypothesisM}
There exist two deterministic positive constants $M$ and $p$
such that
\[\P\left(\sum_z|z|^p\pi_{0z}\leq M^p\right)=1.\]
\end{hypothesisM}

Here, and in the rest of this paper $|\cdot|$ denotes the $l^1$-norm on
$\Z^d$.
Each time hypothesis
(M) is invoked, a range of values will be given for $p$, such as
 $p>1$ or $p\geq2$.

%We will say that (M) is satisfied with $p=\infty$
%when there exists a constant
%$M>0$ such that $\P(\pi_{0z}=0)=1$ if $|z|>M$.

%{\bf One would actually want to assume only
%the necessary and sufficient condition for moving up a level, maybe with
%a necessary and sufficient condition for $E_0(\si_1)<\infty$.
%For the law of large numbers, $\sum_{z\in E}|z|\E(\pi_{0z})<\infty$ should
%be enough, while $\sum_{z\in E}|z|^2\E(\pi_{0z})<\infty$ should be the ultimate
%condition for a central limit theorem.}
%

\section{Invariant measures and the law of large numbers}
\label{invariant}
To prove the law of large numbers
we adopt the  point of view of the particle.
More precisely, we consider the Markov process on $\Omega$ with transition
kernel
\[\pihat(\w,A)=P_0^\w(T_{X_1}\w\in A).\]
To apply the ergodic theorem one needs   an invariant measure
for the above process. Since the state space $\Omega$ is large, there will be
many such measures. The one that is useful for us
will be suitably  ``comparable'' to the
original measure $\P$.
  For integers $n$ define $\sigma$-algebras
${\kS}_n=\sigma(\w_x:x\cdot\uhat\geq n)$.
Our first result is the existence and uniqueness theorem
for the invariant distribution.

\begin{theorem}
\label{Pinfty}
Let $d\geq1$ and consider a product probability measure $\P$ on environments
with a forbidden direction $-\uhat\in\R^d\setminus\{0\}$ as in \eqref{forbidden}.
Assume non-nestling {\rm(N)} in direction $\uhat$,
and the moment hypothesis {\rm(M)} with $p>1$.  Then
there exists a probability measure $\P_\infty$
 on $(\Omega,{\kS})$ that is invariant for the
Markov process with transition kernel $\pihat$ and has these
properties:
\begin{enumerate}
\item[{\rm (a)}]
 $\P_\infty$ is absolutely continuous relative to $\P$ when restricted
to ${\kS}_k$ with  $k\leq0$.  Moreover, $\P=\P_\infty$ on $\kS_1$ and the
two measures are mutually absolutely continuous on $\kS_0$.
\item[{\rm (b)}] The Markov process with kernel $\pihat$ and initial distribution
$\P_\infty$ is ergodic.
\item[{\rm (c)}] There can be at most one $\P_\infty$ that is invariant
under $\pihat$, absolutely continuous relative to $\P$ on
${\kS}_k$ for all $k\leq0$, and satisfies {\rm (b)}. Moreover,
$\P_n(A)=P_0(T_{X_n}\w\in A)$ converges to $\P_\infty(A)$ for any $A\in\kS_k$
for any $k\leq0$.
\end{enumerate}
\end{theorem}

An issue of interest in the RWRE 
literature is the equivalence (mutual absolute continuity)
of  $\P$ and $\P_\infty$. 
To have equivalence  of $\P$ and $\P_\infty$
on all negative half-spaces requires  some sort of ellipticity.
An example of strong enough ellipticity hypotheses that imply
such absolute continuity can be found in  Theorem 2 of \cite{llnrwre}.
On page \pageref{hypoE} below we introduce the ellipticity hypothesis (E)
under which the CLT scenario II of the introduction is valid, as
proved in the companion paper \cite{forbidden-qclt}.  Here we wish to
point out that this  condition (E)  is still too weak to
imply  mutual absolute continuity on all $\kS_k$, as illustrated by
this example. 

\begin{example}
\label{abscont}
Consider a two-dimensional product environment $\P$ whose marginal at $0$ is,
with equal probability, one of $\{p_1,p_2,p_3\}$. Here,
\begin{align*}
p_1=\tfrac12(\delta_{(1,1)}+\delta_{(1,-1)}),~
p_2=\tfrac12(\delta_{(1,1)}+\delta_{(1,0)}),\text{ and }
p_3=\tfrac12(\delta_{(1,-1)}+\delta_{(1,0)}).
\end{align*}
Here again we have finite-size jumps, so that the moment hypothesis (M) is
in force. So is the non-nestling hypothesis (N). As for ellipticity, only the
weak version (hypothesis (E) below) is satisfied.

%Now, $P_0(\si_n=n)=1$ for all $n\geq0$, so
%\eqref{Pinftyformula} gives us  that $\P_\infty=\P_n$ on
%$\kS_{-n}$ for all $n\geq0$.
A computation shows that
\[f_n(\w):=\sum_{x\cdot e_1=-n}P_x^\w(X_n=0)=\frac{d\P_n}{d\P}(\w).\]
Thus, this is an $L^1(\P)$-martingale relative to the filtration
$\{\kS_{-n}\}_{n\geq0}$ and by the martingale convergence theorem it converges
in $L^1(\P)$ to a limit $f$. Then, $\frac{d\P_\infty}{d\P}=f$ and
$\P_\infty=\P_n$ on $\kS_{-n}$ for all $n\geq0$.

Consider the event
$A=\{\w_{(-1,0)}=p_1,\w_{(-1,1)}=p_2,\w_{(-1,-1)}=p_3\}$.
$A$ is $\kS_{-1}$-measurable and
$\P(A)=1/27>0$. But, if $\w\in A$, $n\geq1$, and $x$ is such that
$x\cdot e_1<0$, then $P_x^\w(X_n=0)=0$. This implies
that $\P_\infty(A)=\P_n(A)=0$ for all $n\geq1$.
\end{example}

%Of course, since the walk starting at $0$ stays in $\{x:x\cdot\uhat\geq0\}$,
%the absolute continuity in part (a) above is all one needs.

Once one has a suitable invariant measure one can conclude a law of large
numbers. Define the  {\sl drift} as
 \[D(\w)=E_0^\w(X_1)=\sum_z z\pi_{0z}(\w).\]
By part (a) of Theorem \ref{Pinfty} the moment assumption
(M) is also valid under $\P_\infty$, so there is no problem in
defining   \[v=\E_\infty(D).\]
Here, naturally,  $\E_\infty$ is expectation under
$\P_\infty$ of Theorem {\rm\ref{Pinfty}}.
Define also $\si_1
=\inf\{n:X_n\cdot\uhat\geq X_0\cdot\uhat+1\}$. 
One has the following:

\begin{theorem}
\label{lln}
Let $d\geq1$ and consider a product probability measure $\P$ on environments
with a forbidden direction $-\uhat\in\R^d\setminus\{0\}$ as in \eqref{forbidden}.
Assume non-nestling {\rm(N)} in direction $\uhat$,
and the moment hypothesis {\rm(M)} with $p>1$.
Then the following law of large numbers is satisfied:
\[P_0\left(\lim_{n\to\infty}n^{-1}X_n=v\right)=1.\]
Moreover, $E_0(\si_1)<\infty$, $v=E_0(X_{\si_1})/E_0(\si_1)$, and
\begin{align}
\sup_n |E_0(X_n)-nv|<\infty.
\label{Xn-nv}
\end{align}
\end{theorem}

Theorems \ref{Pinfty} and \ref{lln} correspond essentially to Theorems 2 and 3 of
\cite{llnrwre}. However, we will reprove them since some of
 our assumptions are weaker
than those of \cite{llnrwre}. We start with a lemma.

\begin{lemma}
\label{exponential}
Let $d\geq1$ and consider a $T$-invariant probability measure $\P$ on environments
with a forbidden direction $-\uhat\in\R^d\setminus\{0\}$ as in \eqref{forbidden}.
Assume non-nestling {\rm(N)} in direction $\uhat$,
and the moment hypothesis {\rm(M)} with $p>1$.
Then there exist strictly positive, finite
constants $\bar C_m(M,\delta,p)$, $\hat C_{\bar p}(M,\delta,p)$,
and $\lambda_0(M,\delta,p)$
such that
 for all $x\in\Z^d$, $\lambda\in[0,\lambda_0]$, $n,m\geq0$, 
and $\P$-a.e. $\w$,
\begin{align}
&E_x^\w(|X_m-x|^{\bar p})\leq M^{\bar p}m^{\bar p}\,\text{ for }
1\leq\bar p\leq p,
\label{0}\\
&E_x^\w(e^{-\lambda X_n\cdot\uhat})\leq
e^{-\lambda x\cdot\uhat}(1-\lambda\delta/2)^n,
\label{1}\\
%&P_x^\w(X_n\cdot\,e_1\leq a)\leq e^{\lambda(a-x\cdot e_1)}
%(1-\lambda\delta/2)^n,
%\label{2}\\
%&\sup_{x\in\Z^d}P_x^\w(X_n\cdot\,e_1\leq X_0\cdot e_1-k)\leq
%e^{-\lambda k}(1-\lambda\delta/2)^n,
&P_x^\w(\si_1>n)\leq e^\la(1-\lambda\delta/2)^n,
\label{3}\\
&E_x^\w(\si_1^m)\leq \bar C_m,
\label{4}\\
&E_x^\w(|X_{\si_1}-x|^{\bar p})\leq\hat C_{\bar p}\,\text{ for }
1\leq\bar p<p.
\label{5}
\end{align}
\end{lemma}

\begin{proof}
Note first that $E_x^\w(|X_m-x|^{\bar p})=E_0^{T_x\w}(|X_m|^{\bar p})$,
with a similar
equation for each of the above assertions. Therefore, without loss
of generality, we will assume that $x=0$.
To prove \eqref{0} we write
\[E_0^\w(|X_m|^{\bar p})\leq m^{\bar p-1}
\sum_{j=0}^{m-1}E_0^\w\left(E_0^{T_{X_j}\w}(|X_1|^{\bar p})\right)
\leq M^{\bar p}m^{\bar p}.\]
Concerning \eqref{1}, observe that since $P_0^\w(X_1\cdot\uhat\geq0)=1$
$\P$-a.s., one has
\begin{align*}
E_x^\w(e^{-\lambda X_n\cdot\uhat}|X_{n-1})&=
e^{-\lambda X_{n-1}\cdot\uhat}E_0^{T_{X_{n-1}}\w}(e^{-\lambda X_1\cdot\uhat})\\
&\leq e^{-\lambda X_{n-1}\cdot\uhat}
\left(1-\lambda D(T_{X_{n-1}}\w)\cdot\uhat
+\lambda^pE_0^{T_{X_{n-1}}\w}(|X_1|^p)\right)\\
&\leq e^{-\lambda X_{n-1}\cdot\uhat}(1-\lambda\delta+M^p\lambda^p).
\end{align*}
We have used the fact that $D\cdot\uhat\geq\delta$ and $E_0^\w(|X_1|^p)\leq M^p$,
$\P$-a.s.
Taking now the quenched expectation on both sides and iterating the procedure
proves \eqref{1}. Then \eqref{3} follows immediately. Indeed,
\[P_0^\w(\si_1>n)=P_0^\w(X_n\cdot\uhat<1)\leq
e^\la E_0^\w(e^{-\lambda X_n\cdot\uhat})
\leq e^\la (1-\delta\lambda/2)^n,\]
implying also \eqref{4}. As for \eqref{5}, we have
\begin{align*}
E_0^\w(|X_{\si_1}|^{\bar p})&=\sum_{m\geq1}
E_0^\w(|X_m|^{\bar p},\si_1=m)\\
&\leq\sum_{m\geq1}E_0^\w(|X_m|^p)^{\frac{\bar p}{p}}
P_0^\w(X_{m-1}\cdot\uhat<1)^{\frac{p-\bar p}{p}}\\
&\leq M^{\bar p}e^{\la\frac{p-\bar p}{p}}
\sum_{m\geq1}m^{\bar p}\left((1-\lambda\delta/2)^
{\frac{p-\bar p}{p}}\right)^{m-1}
<\infty.\qedhere
\end{align*}
\end{proof}

Once $P_0(\si_1<\infty)=1$ one can set $\si_0=0$ and
define, by induction, for $k\geq1$:
\begin{align}
\label{sik}
\si_{k+1}=\inf\{n>\si_k:X_n\cdot\uhat\geq X_{\si_k}\cdot\uhat+1\}<\infty,
~P_0\text{-a.s.}
\end{align}
Using ideas from \cite{szlln}, we next show an annealed renewal property.

\begin{proposition}
\label{hittingmoments}
Let $d\geq1$ and consider a product probability measure $\P$ on environments
with a forbidden direction $-\uhat\in\R^d\setminus\{0\}$ as in \eqref{forbidden}.
Assume non-nestling {\rm(N)} in direction $\uhat$,
and the moment hypothesis {\rm(M)} with $p>1$.
Then
\begin{enumerate}
\item[(a)] $E_0(\sigma_k)<\infty$ for all $k\geq0$. Moreover,
\begin{align*}
P_0\left(\lim_{k\to\infty}k^{-1}\sigma_k=E_0(\sigma_1)\right)=1.
\end{align*}
\item[(b)] Fix some integer $K\geq0$ and define
\begin{align*}
\Theta_k=\bigl(\si_{k+1}-\si_k\,;\,
&X_{\si_k+1}-X_{\si_k},X_{\si_k+2}-X_{\si_k},
\dotsc,
X_{\si_{k+1}}-X_{\si_k}\,; \\
&\{\w_x:|x\cdot \uhat-X_{\si_k}\cdot \uhat|\leq K\}\bigr).
\end{align*}
Then $(\Theta_{m(2K+1)+m_0})_{m\geq1}$
is i.i.d. under $P_0$ for any $m_0\geq0$. In particular, if $K=0$, then
$(\Theta_m)_{m\geq0}$ is i.i.d. under $P_0$.
\end{enumerate}
\end{proposition}

\begin{proof}
First, by \eqref{4} we know that $E_0(\si_1)<\infty$.
If one takes $m_0=0$ and $K=0$
in (b), then $(\si_{k+1}-\si_k)_{k\geq0}$ is a sequence of non-negative i.i.d.
random variables under  $P_0$ and the rest of (a) follows. So we will now prove (b).

Fix $k\geq0$ and $L\geq k+1+2K$ and define
\[\cG_k=
\sigma(\si_1,\dotsc,\si_{k+1},X_1,\dotsc,X_{\si_{k+1}},
\{\w_x:x\cdot \uhat\leq X_{\si_k}\cdot\uhat+K\}).\]
Let $F$ be a bounded function on $\Z^d$ paths $\Xhat$ and environments
$\w$ which depends on $\w$ only through $\kS_{-K}$.
Let $H$ be a bounded $\cG_k$-measurable function.

To understand where the independence here and in similar places
elsewhere in the paper comes from, observe that a quenched probability
such as $P^\w_0(X_{\si_k}=x)$ is a function of
$(\w_z: z\cdot\uhat<x\cdot\uhat)$, while any quenched expectation
$E^\w_x[ f(X_0, X_1,X_2,\dotsc)]$ of the process restarted at $x$
 is a function of
$(\w_z: z\cdot\uhat\geq x\cdot\uhat)$.

Since $\si_k<\infty$ $P_0$-a.s., we can write
\begin{align*}
\begin{split}
E_0[F(X_{\si_L+\,\cdot}-X_{\si_L},T_{X_{\si_L}}\w)H]
&=\sum_{x}\E[E_0^\w(H,X_{\si_{L-K}}=x)]\\
&\qquad\times\E[E_x^\w\{F(X_{\si_K+\,\cdot}-X_{\si_K},T_{X_{\si_K}}\w)\}]
\end{split}\\
\begin{split}
&=\sum_{x}E_0[H,X_{\si_{L-K}}=x]\\
&\qquad\times E_0[F(X_{\si_K+\cdot}-X_{\si_K},T_{X_{\si_K}}\w)]
\end{split}\\
\begin{split}
&=E_0[H]E_0[F(X_{\si_K+\cdot}-X_{\si_K},T_{X_{\si_K}}\w)].
\end{split}
\end{align*}
In the first equality we  used
\[X_{\si_k}\cdot\uhat+K<X_{\si_{L-K-1}}\cdot\uhat+1\leq X_{\si_{L-K}}\cdot\uhat
\leq X_{\si_L}\cdot\uhat-K.\]

With $L=m(2K+1)+m_0$ and $k=(m-1)(2K+1)+m_0$ the above shows  that
$\Theta_{m(2K+1)+m_0}$ is independent of $\cG_{(m-1)(2K+1)+m_0}$
and distributed like $\Theta_K$. Since $\Theta_k$ is $\cG_k$-measurable,
the i.i.d. property follows. If $K=0$, then one can start $m$ from $0$,
since $\Theta_K=\Theta_0$.
\end{proof}

The next proposition defines the invariant measure through a limit
and a formula.

\begin{proposition}
\label{Pn}
Let $d\geq1$ and consider a product probability measure $\P$ on environments
with a forbidden direction $-\uhat\in\R^d\setminus\{0\}$ as in \eqref{forbidden}.
Assume non-nestling {\rm(N)} in direction $\uhat$,
and the moment hypothesis {\rm(M)} with $p>1$.
Then there exists a probability measure $\P_\infty$ such that
the sequence $\P_n(A)=P_0(T_{X_n}\w\in A)$
 converges to $\P_\infty(A)$ for all $A\in\kS_{-k}$ and
all $k\geq0$.
Moreover, if $A$ is $\kS_{-k}$-measurable for some $k\geq0$, then
\begin{align}
\label{Pinftyformula}
\P_\infty(A)=\frac{E_0\left(\sum_{m=\si_k}^{\si_{k+1}-1}
\one\{T_{X_m}\w\in A\}\right)}{E_0(\si_1)}.
\end{align}
\end{proposition}

\begin{proof}
Write
\[P_0(T_{X_n}\w\in A)=P_0(\si_{k+1}>n,T_{X_n}\w\in A)+
P_0(\si_{k+1}\leq n,T_{X_n}\w\in A).\]
The first term on the right-hand-side goes to $0$ as $n$ goes to infinity.
As for the second term, we have
\begin{align*}
&P_0(\si_{k+1}\leq n,T_{X_n}\w\in A)=
\sum_{j\geq1}P_0(\si_{k+j}\leq n<\si_{k+j+1},T_{X_n}\w\in A)\\
&\qquad=\sum_{\substack{j\geq1,x,y\\ 0\leq m\leq n}}
P_0(X_m=x,\si_j=m,\si_{k+j}\leq n<\si_{k+j+1},X_n=y,T_y\w\in A)\\
&\qquad=\sum_{\substack{j\geq1,x,y\\ 0\leq m\leq n}}
P_0(X_m=x,\si_j=m)P_0(\si_k\leq n-m<\si_{k+1},X_{n-m}=y-x,T_{y-x}\w\in A)\\
&\qquad=\sum_{\substack{j\geq1\\ 0\leq m\leq n}}
P_0(\si_j=m)P_0(\si_k\leq n-m<\si_{k+1},T_{X_{n-m}}\w\in A)\\
&\qquad=\sum_{m=0}^n P_0(\exists j\geq1:\si_j=n-m)
P_0(\si_k\leq m<\si_{k+1},T_{X_m}\w\in A).
\end{align*}
In the third equality we have used the fact that $T_{y}A$ is
$\kS_{y\cdot \uhat-k}$-measurable, while
\[y\cdot \uhat=X_n\cdot\uhat\geq X_{\si_{k+j}}\cdot\uhat
\geq X_{\si_j}\cdot\uhat+k=x\cdot \uhat+k.\]
Finally, by Proposition \ref{hittingmoments}, $(\si_j-\si_{j-1})_{j\geq1}$ are
independent under $P_0$ and thus, by the renewal theorem
\cite[p. 313]{feller1}, we have that $P_0(\exists j\geq1:\si_j=n-m)$
converges to $E_0(\si_1)^{-1}$. This finishes the proof.
\end{proof}

We are now ready to prove Theorems \ref{Pinfty} and \ref{lln}.

\begin{proof}[Proof of Theorem {\rm\ref{Pinfty}}]
Let $\P_\infty$ be the measure in Proposition \ref{Pn}.
For $k\leq0$ and $A\in\kS_k$, $\pihat(\w,A)$ is also $\kS_k$-measurable.
Then the invariance of $\P_\infty$ follows from passing to the limit
in the equation  $\P_{n+1}=\pihat\P_n$, utilizing  the
 convergence proved in Proposition \ref{Pn}.
Moreover, \eqref{Pinftyformula} shows the absolute
continuity of $\P_\infty$ relative to $\P$ on each $\kS_k$ and the absolute
continuity of $\P$ relative to $\P_\infty$ on $\kS_0$. Now, if $A$ is
$\kS_1$-measurable, then
\[\P_n(A)=P_0(T_{X_n}\w\in A)=\sum_x\int P_0^\w(X_n=x)\one\{T_x\w\in A\}\P(d\w)=
\P(A),\]
since the two integrands above are independent and $\P$ is shift-invariant.
This proves that
$\P=\P_\infty$ on $\kS_1$, finishing the proof of part (a).

Consider a bounded local function $\Psi$ that is
$\si(\w_x:|x\cdot \uhat|\leq K)$-measurable for some $K\geq0$.
Let
$\Phi_k=\sum_{j=\si_k}^{\si_{k+1}-1}\Psi(T_{X_j}\w)$.
Due to Proposition \ref{hittingmoments}, $(\Phi_{m(2K+1)+m_0})_{m\geq1}$
is a sequence of i.i.d. random variables with a finite first moment
under probability $P_0$.
This implies that
\begin{align}
P_0\left(\lim_{m\to\infty}
m^{-1}\sum_{j=0}^{\sigma_{(m+1)(2K+1)}-1}\Psi(T_{X_j}\w)=
(2K+1)E_0(\Phi_K)\right)=1.
\label{ergod}
\end{align}
Note that in \eqref{ergod} the sum started at $0$ instead of $\si_{2K+1}$,
since $P_0(\si_{2K+1}<\infty)=1$ and $\Psi$ is bounded.
Define now $K_n=(m_n+1)(2K+1)$ such that $\si_{K_n}\leq n<\si_{K_{n+1}}$.
Then
\[\abs{n^{-1}\sum_{m=0}^{n-1}\Psi(T_{X_m}\w)-
n^{-1}\sum_{m=0}^{\si_{K_n}-1}\Psi(T_{X_m}\w)}\leq
n^{-1}(\si_{K_{n+1}}-\si_{K_n})\norm{\Psi}_\infty\]
goes to $0$ $P_0$-a.s., since $n^{-1}\si_{K_{n}}$ goes to $1$.
On the other hand, by Proposition \ref{hittingmoments},
$k^{-1}\si_k$ converges to $E_0(\si_1)<\infty$. Hence, $n^{-1}K_n$
converges to $E_0(\si_1)^{-1}$ and $n^{-1}m_n$ converges to
$[(2K+1)E_0(\si_1)]^{-1}$. Thus, one has
\begin{align*}
P_0\left(\lim_{n\to\infty}
n^{-1}\sum_{m=0}^{\si_{K_n}-1}\Psi(T_{X_m}\w)=
E_0(\Phi_K)/E_0(\si_1)\right)=1.
\end{align*}
Therefore, we have $\P$-a.s.
\begin{align}
P_0^\w\left(\lim_{n\to\infty}
n^{-1}\sum_{m=0}^{n-1}\Psi(T_{X_m}\w)=E_0(\Phi_K)/E_0(\si_1)\right)=1.
\label{ergod1}
\end{align}
Since the above quantity is $\kS_{-K}$-measurable,
the same holds $\P_\infty$-a.s.
Moreover, the constant limit cannot be anything other than $\E_\infty(\Psi)$.
By bounded convergence, one then has
\begin{align*}
\lim_{n\to\infty}
n^{-1}\sum_{m=0}^{n-1}E_0^\w[\Psi(T_{X_m}\w)]=\E_\infty(\Psi),~\P_\infty\text{-a.s.}
%\label{ergod1}
\end{align*}
Now approximate a general $\Psi\in L^1(\P_\infty)$ by a bounded
local $\Psi$  in the $L^1(\P_\infty)$-sense.
The above convergence
is then  true in the $L^1(\P_\infty)$-sense  for any
$\Psi\in L^1(\P_\infty)$.  The fact that all these limits
are constants suffices for ergodicity. This follows from
Section IV.2 in Rosenblatt's monograph \cite{rosenblatt},
see especially  Corollary 2,
the definition at the  top of p.~94, and Corollary 5.
Part (b) is proved.
%in fact: once we have (\ref{birkhoff}), we know it's true for all bounded functions
%by L^1 approximations. then taking \Psi(\w)=P_0^\w(B) with B time-shift-invariant
%implies that \Psi(T_{X_n}\w) is a martingale an converges to \one_B. on the other hand,
%(\ref{birkhoff}) implies that the limit is a constant P_0(B). thus P_0(B) is 0 or 1.

Once one has ergodicity and absolute continuity on half-spaces
uniqueness follows. Indeed, recalling that \eqref{ergod1} holds $\P$-a.s. and
using bounded convergence once again implies that for any measurable bounded $\Psi$
\[\lim_{n\to\infty} n^{-1}\sum_{m=0}^{n-1}E_0(\Psi(T_{X_m}\w))
=\E_\infty(\Psi),\]
determining $\P_\infty$ uniquely.
The rest of part (c) is true by construction of $\P_\infty$.
\end{proof}

\begin{proof}[Proof of Theorem {\rm\ref{lln}}]
We first prove a law of large numbers with limiting velocity
\[
\bar v=
\frac{E_0(X_{\si_1})}{E_0(\si_1)}
\]
and then the approximation   \eqref{Xn-nv} with $\bar v$ instead of $v$.
Lastly we identify $\bar v$ with $v=\E_\infty(D)$.

Recall that $E_0(|X_{\si_1}|)<\infty$.
Using $K=0$ and $m_0=0$
in part (b) of Proposition \ref{hittingmoments} one can
see that $(X_{\si_k}-X_{\si_{k-1}},\si_k-\si_{k-1})_{k\geq1}$ is a
sequence of i.i.d.
random variables, under $P_0$. Therefore, we have
\[P_0\left(\lim_{k\to\infty}k^{-1}(X_{\si_k},\si_k)=
(E_0(X_{\si_1}),E_0(\si_1))\right)=1.\]
Now let $K_n=K(n)=\max\{j:\sigma_j\leq n\}$. We will sometimes use the notation
$K(n)$ to avoid subscripts on subscripts. Clearly,  $n^{-1}K_n$ converges $P_0$-a.s.
to $E_0(\si_1)^{-1}$.

Write
\begin{align}
\frac{X_n}n=\frac{X_n-X_{\si_{K(n)+1}}}{n}+
\frac{X_{\si_{K(n)+1}}}{K_n+1}\frac{K_n+1}n.
\label{Xn/n}
\end{align}
Using \eqref{5} with some $\bar p\in(1,p)$, we have that for any $\e>0$
\[P_0(|X_n-X_{\si_{K(n)+1}}|>\e n)\leq
\e^{-\bar p}n^{-\bar p}E_0(|X_n-X_{\si_{K(n)+1}}|^{\bar p})
\leq \e^{-\bar p}\hat C_{\bar p} n^{-\bar p}\]
is the general term of a convergent series. By Borel-Cantelli's Lemma, the first term
on the right-hand-side of \eqref{Xn/n} goes to $0$ $P_0$-a.s.
Therefore,
\[P_0\left(\lim_{n\to\infty}n^{-1}X_n=\bar v\right)=1.\]
This proves the law of large numbers with $\bar v$ instead of $v$.
%Note first that by the ergodic theorem one has for $\P_\infty$-a.e. $\w$
%\begin{align*}
%P_0^\w\left(
%\lim_{n\to\infty}n^{-1}\sum_{k=0}^{n-1}D(T_{X_n}\w)
%=\E_\infty(D)\right)=1.
%\end{align*}
%Due to the absolute continuity of $\P$ relative to $\P_\infty$ on $\kS_0$
%the above is also true $\P$-a.s.
%But then the moment hypothesis (M) with $p\geq2$ implies that
%\[M_n=X_n-\sum_{k=0}^{n-1}D(T_{X_n}\w)\] is an
%$L^2(P_0^\w)$-martingale for $\P$-a.e. $\w$.
%Letting $n_j=j^r$, for $r>1$, and using Doob's inequality one has
%\[P_0^\w\left(\max_{k\leq n_j}|M_k|\geq\e n_j\right)
%\leq4(\e n_j)^{-2} E_0^\w(|M_{n_j}|^2)=\Ord(j^{-r}).\]
%By Borel-Cantelli's Lemma $n_j^{-1}\max_{k\leq n_j}|M_k|$ converges to $0$
%$P_0^\w$-a.s. But if $n_j\leq n\leq n_{j+1}$, then
%\[n^{-1}\max_{k\leq n}|M_k|\leq\frac{n_{j+1}}{n_j}\
%n_{j+1}^{-1}\max_{k\leq n_{j+1}}|M_k|\]
%and thus we know that $n^{-1}M_n$ converges to $0$ $P_0$-a.s.
%This proves the law of large numbers.
Since
\begin{align*}
\{K_n=k\}=\left\{\sum_{j=1}^k (\sigma_j-\sigma_{j-1})
\leq n < \sum_{j=1}^{k+1} (\sigma_j-\sigma_{j-1})\right\},
%\label{Kn-def}
\end{align*}
$K_n+1$ is a stopping time  with respect to the filtration of the
process  $\{(X_{\sigma_{k}}-X_{\sigma_{k-1}},
\sigma_{k}-\sigma_{k-1}): k\geq 1\}$.  Moreover, $K_n+1$ is bounded, and
both $\sigma_1$ and
$X_{\sigma_1}$ are integrable by Lemma \ref{exponential}.
By Wald's  identity,
\[E_0\left(\sum_{j=1}^{K_n+1}(X_{\sigma_j}-X_{\sigma_{j-1}})\right)
= E_0(1+K_n) E_0(X_{\sigma_1})\]
and
\[E_0(\sigma_{K_n+1})=
E_0\left(\sum_{j=1}^{K_n+1} ({\sigma_j}-{\sigma_{j-1}})\right)
 = E_0(1+K_n) E_0({\sigma_1}).\]
From
$X_n = \sum_{j=1}^{K_n+1} (X_{\sigma_j}-X_{\sigma_{j-1}})
-(X_{\sigma_{K_n+1}}-X_n)$
we get
\begin{align*}
E_0(X_n) &= E_0(1+K_n) E_0(X_{\sigma_1})
-E_0 (X_{\sigma_{K_n+1}}-X_n)\\
&=
n\bar v + \bar v\left[E_0(1+K_n) E_0(\sigma_1)-n\right]
-E_0(X_{\sigma_{K_n+1}}-X_n)\\
&=
n\bar v + \bar v E_0(\sigma_{K_n+1}-n) -E_0 (X_{\sigma_{K_n+1}}-X_n).
\end{align*}
The last two terms above are bounded by a constant.
To see it for the last one, note that
\begin{align*}
E^\omega_{0}(X_{\sigma_{K_n+1}}-X_n)=
E^\omega_0\left[E^\omega_{X_n}(X_{\sigma_{1}}-X_0)\right]
= E^\omega_0\left[E^{T_{X_n}\omega}_{0}(X_{\sigma_{1}})\right]
\end{align*}
which is bounded by a constant $\P$-almost surely by \eqref{5}.
The same argument works for $E_0(\sigma_{K_n+1}-n)$ via \eqref{4}.
This proves \eqref{Xn-nv} with $\bar v$ instead of $v$.
One consequence of this is that $n^{-1}E_0(X_n)$ converges to $\bar v$, as $n$
goes to infinity.

On the other hand, by the ergodic theorem one has for $\P_\infty$-a.e. $\w$
\begin{align*}
P_0^\w\left(
\lim_{n\to\infty}n^{-1}\sum_{m=0}^{n-1}D(T_{X_m}\w)
=v=\E_\infty(D)\right)=1.
\end{align*}
Due to the absolute continuity of $\P$ relative to $\P_\infty$ on $\kS_0$
the above is also true $\P$-a.s. Due to the moment hypothesis (M), the drift
$D$ is bounded in norm by $M$, for $\P$-a.e. $\w$. Therefore, by bounded convergence,
we conclude that $n^{-1}E_0(X_n)=n^{-1}E_0(\sum_{m=0}^{n-1}D(T_{X_m}\w))$ converges
to $v$, $P_0$-a.s. Thus $\bar v=v$
and the proof of Theorem \ref{lln} is complete.
\end{proof}

\section{The annealed invariance principle}
\label{inv1}
Let us start with more notation and definitions.
We write $\Gamma^t$ for the transpose of a vector or matrix
$\Gamma$. An element of $\R^d$ is regarded as a $d\times 1$ matrix,
or column vector.

For a symmetric, non-negative definite $d\times d$ matrix $\Gamma$,
a Brownian motion with diffusion matrix $\Gamma$ is the $\R^d$-valued
process $\{W(t):t\geq0\}$
such that $W(0)=0$, $W$ has continuous paths, independent increments,
and for $s<t$ the $d$-vector $W(t)-W(s)$ has Gaussian distribution
with mean zero and covariance matrix
$(t-s)\Gamma$.
%If the rank of $\Gamma$ is $m$,
%one can produce such a process by
%finding a $d\times m$ matrix $\Lambda$ such that
%$\Gamma=\Lambda\Lambda^t$, and defining $W(t)=\Lambda B(t)$
%where $B$ is an  $m$-dimensional standard Brownian motion.
The matrix $\Gamma$  is {\it degenerate} in direction $\xi\in\R^d$
if $\xi^t\Gamma \xi=0$.   Equivalently,
$\xi\cdot W(t)=0$ almost surely.

For $t\geq0$  let
\begin{align*}
B_n(t)=\frac{X_{[nt]}-[nt]v}{\sqrt{n}}.
\end{align*}
Here  $[x]=\max\{n\in\Z: n\leq x\}$ for $x\in\R$.
$D_{\R^d}([0,\infty))$ denotes the space of $\R^d$-valued cadlag paths 
endowed with the usual
Skorohod topology, as developed for instance in  \cite{EK}.
%For a closed interval $I\subset[0,\infty)$ denote by $D_{\R^d}(I)$
%the space of right continuous functions on $I$ taking values in $\R^d$
%and having left limits.
%The space $D_{\R^d}(I)$ is endowed with the usual
%Skorohod topology \cite{EK}.
The next theorem is the  annealed invariance principle. 

%\begin{theorem}
%\label{a-clt} Let $d\geq1$ and consider a product probability
%measure $\P$ on environments with a forbidden direction
%$-\uhat\in\R^d\setminus\{0\}$ as in \eqref{forbidden}. Assume that
%$E_0(\si_1^2)<\infty$ and $E_0(\sup_{j\leq\si_1}|X_j|^2)<\infty$.
%Then the distribution of the process $B_n(t)$ under $P_0$
%converges weakly to the distribution of a Brownian motion with
%diffusion matrix
%\begin{align*}
%%\label{diffmat}
%{\mathfrak D}=\frac{E_0[(X_{\si_1}-v\si_1)(X_{\si_1}-v\si_1)^t]}{E_0[\si_1]}.
%\end{align*}
%The matrix $\mathfrak D$ is degenerate in direction $u$ if, and only if,
%$u$ is orthogonal to
%the vector space spanned by $\{x-y:\E(\pi_{0x})\E(\pi_{0y})>0\}$.
%\end{theorem}

\begin{theorem}
\label{a-clt}
Let $d\geq1$ and consider a product probability measure $\P$ on environments
with a forbidden direction $-\uhat\in\R^d\setminus\{0\}$ as in \eqref{forbidden}.
Assume non-nestling {\rm(N)} in direction $\uhat$,
and the moment hypothesis {\rm(M)} with $p>2$.
Then the distribution of the process $B_n(t)$ under $P_0$
converges weakly on the space $D_{\R^d}([0,\infty))$
 to the distribution of a Brownian motion with
diffusion matrix
\begin{align*}
%\label{diffmat}
{\mathfrak D}=\frac{E_0[(X_{\si_1}-v\si_1)(X_{\si_1}-v\si_1)^t]}{E_0[\si_1]}.
\end{align*}
The matrix $\mathfrak D$ is degenerate in direction $u$ if, and only if,
$u$ is orthogonal to
the vector space spanned by $\{x-y:\E(\pi_{0x})\E(\pi_{0y})>0\}$.
\end{theorem}

Note that degeneracy in directions that are orthogonal to all $x-y$,
where $x$ and $y$ range over
admissible jumps, cannot be avoided. This can be seen from
the simple example
of a homogeneous random walk that chooses
 with equal probability between two jumps $a$ and  $b$.
The diffusion matrix is then   $\tfrac14(a-b)(a-b)^t$.

\begin{proof}[Proof of Theorem {\rm\ref{a-clt}}]
Using Lemma \ref{exponential} one can check that
non-nestling (N) in direction $\uhat$ and the moment hypothesis (M) with $p>2$
imply that $E_0(\si_1^2)<\infty$ and $E_0(\sup_{j\leq\si_1}|X_j|^2)<\infty$.
Then, since $X_{\si_k}-v\si_k$ is a sum of square-integrable i.i.d. random variables under
$P_0$, the process $n^{-1/2}( X_{\si_{[nt]}}-v\si_{[nt]})$ converges to a 
Brownian motion with diffusion matrix
$E_0[(X_{\si_1}-v\si_1)(X_{\si_1}-v\si_1)^t]$. 

Next the invariance principle is transferred to the process
$n^{-1/2}(X_{\si_{K([nt])}}-v\si_{K([nt])})$, where as before  
$K_n=K(n)=\max\{j:\sigma_j\leq n\}$.
Now the  limiting diffusion matrix  is 
$\mathfrak D$. We omit this step as it can be found in the proof of 
 Theorem 4.1 of \cite{slowdown+clt}. Note that the
regeneration times $\tau_k$ in \cite{slowdown+clt}
 are exactly our $\si_k$, and that the event
$D=\infty$ in \cite{slowdown+clt} is of full $P_0$-measure in our case.

The last step passes the invariance principle to the process $B_n(t)$,
via the following estimation. 
\begin{align*}
&P_0\left[\sup_{t\leq1}|B_n(t)-n^{-1/2}(X_{\si_{K([nt])}}-v\si_{K([nt])})|
\geq\e\right]\\
&\quad\leq
P_0\left[\sup_{k\leq K_n}\sup_{\si_k\leq j<\si_{k+1}}|X_j-X_{\si_k}|
\geq\e\sqrt n/2\right]+
P_0\left[\sup_{k\leq K_n}(\si_{k+1}-\si_k)\geq\e\sqrt n/2|v|\right]\\
&\quad\leq (n+1)P_0\left[\sup_{0\leq j<\si_1}|X_j|\geq\e\sqrt n/2\right]+
     (n+1)P_0[\si_1\geq\e\sqrt n/2|v|]\\
&\quad\leq\frac{4(n+1)}{n\e^2}\left\{E_0\left[\sup_{0\leq j<\si_1}|X_j|^2,
\sup_{0\leq j<\si_1}|X_j|\geq\e\sqrt n/2\right]+
|v|^2 E_0[\si_1^2,\si_1\geq\e\sqrt n/2]\right\}\\
&\quad\mathop{\longrightarrow}_{n\to\infty}0.
\end{align*}
We used the stationarity of
$(\si_{k+1}-\si_k,X_{\si_k},\dotsc,X_{\si_{k+1}-1})_{k\geq0}$
and
$K_n\leq n$ which follows from $\si_n\geq n$.
In the last step we used 
 $E_0(\si_1^2)<\infty$ and $E_0(\sup_{0\leq j<\si_1}|X_j|^2)<\infty$.
This proves the weak convergence of the process $B_n(t)$ under $P_0$. 

%When Alain switches from inv princ with time [nt] to that with time k_[nt]
% he uses:
%Let X_n converge weakly to X on D-space, assume X has continuous
%paths. Let S_n be a sequence of nondecreasing cadlag processes
%such that for each fixed t,  S_n(t)\to at in probability. Then
%X_n(S_n(t)) converges weakly to X(at) on D-space.
% in that particular case:
%P_0(|\Si_{k_{[nt]}}-\Si_{cnt}|>\e\sqrt n)\leq
%             P_0(|k_{[nt]}-cnt|>\del n)+
%             P_0(\max_{|i|\leq\del n}|\sum_{j=0}^i Z_j|>\e\sqrt n)
%where Z_j are i.i.d.
%by (4.5), the first term goes to 0 for any fixed \del>0.
%by the regular invariance principle, the second goes to 0 after taking n
%to infinity, and then \del to 0.
%this shows the one time marginals converge. for the others, i would guess
%one can break it down in a similar way.

Next we address the  degeneracy of the diffusion matrix. We will
first show that this matrix is degenerate in direction $u$ if, and
only if, $u$ is orthogonal to the vector space spanned by
$\{x-v:\E(\pi_{0x})>0\}$. To this end, let $u\neq0$ be such that
$P_0(X_{\si_1}\cdot u=\si_1 v\cdot u)=1$.

If $x$ is such that $x\cdot\uhat>0$ and $\E_0(\pi_{0x})>0$, then
\[P_0(X_{\si_1}=x,\si_1=1)\geq\E(\pi_{0x})>0\]
implies that $x\cdot u=v\cdot u$ and $u$ is perpendicular to
$x-v$.

If $z\not=0$ is such that $z\cdot\uhat=0$ and $\E(\pi_{0z})>0$,
then choose $x$ such that $x\cdot\uhat>0$ and $\E(\pi_{0x})>0$.
This is always possible, due to non-nestling (N). Then
\[P_0(X_{\si_1}=nz+x,\si_1=n+1)\geq \E(\pi_{0z})^n\E(\pi_{0x})>0\]
implies that $nz\cdot u+x\cdot u=(n+1)v\cdot u$, for all $n\geq0$.
This gives $z\cdot u=x\cdot u=v\cdot u$ and $u$ is perpendicular
to $z-v$.

If $\E(\pi_{00})>0$, then by non-nestling (N) there exists a vector
$z\neq0$ such that $\E(\pi_{00}\pi_{0z})>0$. If $z\cdot\uhat>0$,
then
\[P_0(X_{\si_1}=z,\si_1=n+1)\geq \E(\pi_{00}^n\pi_{0z})>0\]
and $z\cdot u=(n+1)v\cdot u$ implies that $v\cdot u=0$.  On the
other hand, if $z\cdot\uhat=0$, then choose $x$ such that
$x\cdot\uhat>0$ and $\E(\pi_{0x})>0$. In this case, we have
\[P_0(X_{\si_1}=z+x,\si_1=n+2)\geq \E(\pi_{00}^n\pi_{0z})\E(\pi_{0x})>0,\]
which again implies $(z+x)\cdot u=(n+2)v\cdot u$ and $v\cdot u=0$.

The above discussion shows that if the diffusion matrix is
degenerate in direction $u$, then $u$ is orthogonal to the vector
space spanned by $\{x-v:\E(\pi_{0x})>0\}$. Inversely, if $u$ is
orthogonal to $x-v$ for each $x$ such that $\E(\pi_{0x})>0$, then
\[E_0[|u\cdot (X_{\si_1}-\si_1v)|]
=\sum_{n\geq1}E_0\left[\one\{\si_1=n\} \abs{\sum_{k=0}^{n-1}
u\cdot((X_{k+1}-X_k)-v)}\right]=0\] and the diffusion matrix is
degenerate in direction $u$.

To finish, we observe that the span of $\{x-v:\E(\pi_{0x})>0\}$ is
the same as that of $\{x-y:\E(\pi_{0x})\E(\pi_{0y})>0\}$, which is
the space that appears in the theorem. One part comes simply  from
$x-y=(x-v)-(y-v)$. Conversely,
\[x-v=\lim_{n\to\infty}n^{-1}\sum_{k=0}^{n-1}(x-(X_{k+1}-X_k))\]
shows that $x-v$ lies in the span of $x-y$ for admissible steps
$y$.
%(Recall that vector subspaces in finite dimensions are closed.)
\end{proof}

\section{The quenched invariance principle}
\label{inv2}
%By the bounded convergence theorem one sees that a quenched invariance principle
%would imply an annealed one.
Two types of centering of the quenched process will be considered: the (random) quenched mean
$E_0^\w(X_n)$ and its asymptotic limit $nv$.
For $\w\in\Omega$,  define two scaled processes
\begin{align*}
B_n(t)=\frac{X_{[nt]}-[nt]v}{\sqrt{n}}\quad\text{ and }\quad
\Btil_n(t)=\frac{X_{[nt]}-E_0^\w(X_{[nt]})}{\sqrt{n}}.
\end{align*}
Let $Q_n^\w$, respectively $\widetilde Q_n^\w$, denote the distribution of $B_n$,
respectively $\Btil_n$, induced by $P_0^\w$ on the Borel sets of
$D_{\R^d}([0,\infty))$.   The key to distinguishing between the two 
quenched scenarios I and II described in the Introduction is the following
ellipticity hypothesis.

\newtheorem*{hypothesisE}{\sc Hypothesis (E)}
\begin{hypothesisE}
\label{hypoE}
One has %$\P\left(\sup_z\pi_{0z}<1\right)>0$ and
\begin{align}
\label{badcase}
\P(\,\forall z\not=0:\pi_{00}+\pi_{0z}<1)>0.
\end{align}
Moreover, the walk is not supported by any one-dimensional subspace.
More precisely, if $\cJ=\{y:\E(\pi_{0y})>0\}$ is the set
of all points that are accessible from $0$ with one jump,
then $\cJ$ is not contained in any subspace of the kind
$\R u=\{ su: s\in\R\}$ for any $u\in\R^d$.
In particular, this rules out the case $d=1$.
\end{hypothesisE}

It turns out that if we add assumption (E) to the earlier hypotheses 
\eqref{forbidden}, (N), and (M) with $p>2$, 
then both  laws $Q_n^\w$ and $\widetilde Q_n^\w$
converge to a Brownian motion with the diffusion matrix $\kD$ of Theorem \ref{a-clt},
for $\P$-almost every $\w$.  This is proved in the companion paper 
 \cite{forbidden-qclt}. 

What we show in the remainder of the present paper is that 
(E) is also {\sl necessary}  in order for the processes $B_n$ 
to satisfy  the quenched invariance principle. 
So we have an almost complete characterization of the central limit
behavior of the non-nestling walk with a forbidden direction. The 
only shortcoming is that some of our  moment
 hypotheses are 
 almost certainly not best possible. 

To be more precise, we show that if (E) fails while the other 
hypotheses are in force,  then scenario I from the Introduction 
is realized.  And additionally,  if a degenerate situation is ruled 
out, the distributions $\{Q_n^\w\}_{n\geq 1}$ are not tight 
for a.e. $\w$, 
so the processes $B_n$
cannot satisfy a quenched central limit theorem. 
The proof is done separately for  the two possible ways in which
 hypothesis (E) can fail:

(i) The walk can be supported on a one-dimensional subspace.
Such a walk can be converted to the case $d=1$, which is treated in
Theorem \ref{qclt-1d}. 

(ii) We can have $d\geq 2$ but condition
\eqref{badcase} fails.  This case is treated 
 in Theorem \ref{hypE}.

%In both examples, the quenched invariance principle
%does not hold with centering $nv$.
The qualitative insight we take from this 
 is that the quenched CLT centered at $nv$  demands  enough random
noise in the selection of the path
of the walker under a fixed environment. Only then 
can the   walker experience a sufficiently rich
 sample of environments so that the fluctuations of the quenched
mean do not interfere at the diffusive scale. 
%Otherwise, as it is the case in both Theorems,
%the quenched central limit theorem centered at $nv$ cannot hold.
%However,  there is still a quenched central limit theorem
%centered at $E_0^\w(X_n)$.

%Then, the following theorem and its corollary are proved in
%\cite{forbidden-qclt}.
%
%\newtheorem*{theoremQIP}{\sc Theorem Q-IP}
%\begin{theoremQIP}
%Let $d\geq2$ and consider a product probability measure $\P$ on environments
%with a forbidden direction $-\uhat\in\Q^d\setminus\{0\}$ as in \eqref{forbidden}.
%Assume non-nestling {\rm(N)} in direction $\uhat$,
%the moment hypothesis {\rm(M)} with $p>2$, and ellipticity
%{\rm(E)}.  Then for $\P$-a.e. $\w$ the
%distributions $Q_n^\w$ converge weakly to the distribution of a Brownian
%motion with diffusion matrix $\mathfrak D$ as in Theorem \ref{a-clt}.
%
%Moreover, for $\P$-almost every $\w$,
%$n^{-1/2}\max_{k\leq n}|E_0^\w(X_k)-kv|$ converges
%to $0$ as $n$ goes to infinity and, therefore, the
%same invariance principle holds for $\widetilde Q_n^\w$.
%\end{theoremQIP}
%
%
%\newtheorem*{corollaryQIP}{\sc Corollary Q-IP}
%\begin{corollaryQIP}
%Let $d\geq2$ and consider a product
%probability measure $\P$ on environments
%with a forbidden direction $-\uhat\in\R^d\setminus\{0\}$ as in
%\eqref{forbidden}.  Assume the step distribution is finitely
%supported, in other words that the set $\cJ=\{y:\E(\pi_{0y})>0\}$ is finite.
%Assume non-nestling {\rm(N)} in direction $\uhat$  and
%ellipticity {\rm(E)}.
%Then all the conclusions of Theorem Q-IP  hold.
%\end{corollaryQIP}

\subsection{The one-dimensional case}
When $d=1$ we take $\uhat=1$, and so our 
 random walk is totally asymmetric in the sense that
$X_n\leq X_{n+1}$ holds almost surely. 
Define the events $V_i=\{\exists n\geq0:X_n=i\}$ and
the coefficients
\[a_i(\w)=(1-\pi_{ii}(\w))^{-1}(P_0^\w(V_i)-P_1^\w(V_i)).\]
Let $g(\w)=D(\w)-v$ and define
\begin{align}
\label{Delta}
\Delta(\w) = \sum_{i\geq 0} a_i(\w) g(T_i\w).
\end{align}
Two observations that will help seeing how various formulas emerge:
in the one-dimensional case
\[\frac{d\P_\infty}{d\P}=\frac{(1-\pi_{00})^{-1}}{\E[(1-\pi_{00})^{-1}]}\,,\]
 and also
\[\E[|E_0^\w(X_{\si_1})-v E_0^\w(\si_1)|^2]=\E[g^2/(1-\pi_{00})^2].\]
The following constants will be the diffusion coefficients in our next
theorem.
\begin{align}
%=v\E[|E_0^\w(X_{\si_1})-v E_0^\w(\si_1)|^2]/E_{00}(L)
\kappa_m^2&=v^{-1}
\E[|E_0^\w(X_{\si_1})-v E_0^\w(\si_1)|^2]\sum_{i\geq0}\E\left[\Big|
P_0^\w(V_i)-\sum_{k=0}^i P_k^\w(V_i)\E_\infty\pi_{0k}
\Big|^2\right],
\label{kappa_m}\\
%\kappa_q^2&=\bigl(\E[(1-\pi_{00})^{-1}]\bigr)^{-1}
%\E\left[(1-\pi_{00}(\w))^{-1}E_0^\w\Big(\,\Big|X_1-v-
%\sum_{z=0}^{X_1-1}\Delta(T_z\w)\Big|^2\,\Big)\right].
%%=E_0[|X_{\si_1}-v\si_1|^2]/E_0(\si_1)-\kappa_m^2.
%\kappa_q^2&=\bigl(\E[(1-\pi_{00})^{-1}]\bigr)^{-1}
%\E\left[\frac{E_0^\w\Big(\,\Big|X_1-v-
%\sum_{z=0}^{X_1-1}\Delta(T_z\w)\Big|^2\,\Big)}{1-\pi_{00}(\w)}\right].
\kappa_q^2&=\E_\infty E_0^\w\Big(\,\Big|X_1-v-
\sum_{z=0}^{X_1-1}\Delta(T_z\w)\Big|^2\,\Big).
\label{kappa_q}
\end{align}

We first state some properties of the diffusion constants.

\begin{proposition}
\label{km+kq}
Let $d=1$ and consider a product probability
measure $\P$ on environments with forbidden direction $-\uhat=-1$
as in \eqref{forbidden}.  Assume non-nestling {\rm(N)} in direction
$\uhat$, and the moment hypothesis {\rm(M)} with $p>2$. Then
\begin{enumerate}
\item $\kappa_m>0$ if, and only if,
\begin{align}
\label{non-deg1}
\P(D=v)<1.
\end{align}
Moreover, when $\kappa_m=0$, one has $E_0^\w(X_n)=nv$ and
$Q_n^\w=\Qtil_n^\w$, $\P$-a.s. \item $\kappa_q>0$ if, and
only if,
\begin{align}
\label{non-deg2}
\P(\sup_z\pi_{0z}<1)>0.
\end{align}
\item One has \[\kappa_m^2+\kappa_q^2=E_0[|X_{\si_1}-v\si_1|^2]/E_0(\si_1),\]
which is precisely $\mathfrak D$ from Theorem {\rm\ref{a-clt}}.
\item We have this alternative representation: 
\[\kappa_m^2=\frac{v}{E_{0,0}(L)}\E[|E_0^\w(X_{\si_1}-v\si_1)|^2]\]
where $E_{0,0}(L)$ is the expected first common point
\[L=\inf\{\ell>0: \exists n,m: X_n=\bar X_m=\ell\}\]
of two walks that start at the origin and run independently 
in a common environment. 
\end{enumerate}
\end{proposition}

We do not have an argument for the decomposition in (c) directly
from the defining formulas. Instead, we justify it by a limit. 
In order not to interrupt the flow of ideas we will prove this proposition after the statement
and proof of the next theorem. But before we proceed with the theorem,
let us give two examples where one can compute the
above diffusion coefficients.

\begin{example}
Consider a one-dimensional product environment $\P$ with marginal satisfying:
\[\P(\pi_{01}+\pi_{02}=1)=1.\]
For convenience, we will use $p_x=\pi_{x,x+1}$ and $q_x=\pi_{x,x+2}=1-p_x$.
Then,
\[P_0^\w(V_i)-P_1^\w(V_i)=-q_0(P_1^\w(V_i)-P_2^\w(V_i))=\cdots=(-1)^iq_0\cdots q_{i-1}\]
and $v=\E D=1+\E q_0$. Hence
\begin{align*}
\kappa_m^2&=(1+\E q_0)^{-1}\VVar(E_0^\w X_1)\sum_{i\geq0}\E[|
(P_0^\w(V_i)-P_1^\w(V_i))+(P_1^\w(V_i)-P_2^\w(V_i))\E q_0|^2]\\
&=(1+\E q_0)^{-1}\VVar(q_0)\big(1+\sum_{i\geq1}\E[|
(q_0-\E q_0)q_1\cdots q_{i-1}|^2]\big)\\
&=\VVar(q_0)\frac{1-\E q_0}{1-\E(q_0^2)}
\end{align*}
while
\[\kappa_q^2=E_0[|X_1-v|^2]-\kappa_m^2.\]
\end{example}

The next example is used in the proof of Theorem \ref{hypE}.

\begin{example}
Consider a one-dimensional product environment $\P$ with marginal satisfying 
\[\P(\pi_{00}+\pi_{01}=1)=1.\]
Then $v=1/E_0(\si_1)=(E_0[\pi_{01}^{-1}])^{-1}$ and $a_0(\w)=\pi_{01}^{-1}(\w)$ while
$a_i(\w)=0$ for $i\geq1$. This gives
\begin{align}
\kappa_m^2&%=\frac{\E[|E_0^\w(X_{\si_1}-v\si_1)|^2]}{E_0(\si_1)}
=(\E[\pi_{01}^{-1}])^{-3}\VVar(\pi_{01}^{-1})
=\frac{\VVar(E_0^\w[v\si_1])}{E_0(\si_1)},
\label{qmean-cov-la-1}\\
\kappa_q^2&=(\E[\pi_{01}^{-1}])^{-3}
\E\bigl[\frac{\pi_{00}}{\pi_{01}^2}\bigr]
=\frac{\E[\Var^\w_0(v\si_1)]}{E_0(\si_1)}.
\label{ell-var}
\end{align}
\label{1d-nn-example}
\end{example}

The following theorem treats the one-dimensional situation.
The moment needed is defined by 
\be
p_0=\frac{19}6+\frac{\sqrt{139}}3
\cos(\frac13\arccos(\frac{1504}{139^{3/2}})).
\label{def-p0}
\ee  
This  is approximately 7.06025. 

\begin{theorem}
\label{qclt-1d} Let $d=1$ and consider a product probability
measure $\P$ on environments with forbidden direction $-\uhat=-1$
as in \eqref{forbidden}. Assume non-nestling {\rm(N)} in direction
$\uhat$ and the moment hypothesis {\rm(M)} with 
$p>p_0$.
%=\frac{19}6+\frac{\sqrt{139}}3\cos(\frac13\arccos(\frac{1504}{139^{3/2}}))$.
%%%$p>\tfrac12(5+\sqrt{17})$.  %%%$p\geq8$. 
Then
\begin{enumerate}
\item For $\P$-almost every  $\w$ the
distributions $\Qtil_n^\w$ converge weakly to the distribution of
the process $\{\kappa_q B(t) : t\geq 0\}$
where $B(t)$ denotes standard one-dimensional Brownian motion.
\item Under $\P$, the scaled quenched mean process
$n^{-1/2}\{E^\w_0(X_{[nt]})-ntv\}$ converges in distribution
to $\{\kappa_m B(t):t\geq0\}$.
\item Assume the non-degeneracy condition \eqref{non-deg1} holds.
Then  for $\P$-almost every  $\w$, the sequence
of distributions $\{Q^\w_n:n\geq 1\}$ is not tight.
\end{enumerate}
\end{theorem}

%\begin{remark}
%Of course, when $\P(D=v)=1$, $\Q^\w_n=\Qtil^\w_n$.
%\end{remark}

\begin{remark}
The crude moment bounds we use in the proof require the strong assumption
$p>p_0$.
%$p>\tfrac12(5+\sqrt{17})$. 
We believe this 
can be improved, but do not venture into this.
\end{remark}

\begin{proof}%[Proof of Theorem {\rm\ref{qclt-1d}}]
Given the moment $p$ of the hypothesis of the theorem, 
henceforth we write 
\[m=2\bar p-2\] 
where $\bar p$ is some value that satisfies $2\leq \bar p<p$. The lower bound
on $\bar p$ ensures  $m\geq 2$, which is convenient so that powers $m/2$
used below do not go below 1.  Also, $m>1$ is needed for applying
Theorem 3.2 from \cite{burk} in the next proof. 

\begin{lemma}
\label{Delta-moments}
$\Delta$ is $\kS_0$-measurable, $\E(\,\vert\Delta\rvert^m\,)<\infty$,
 and 
$\E(\Delta)=0$.
\end{lemma}
%
%{\bf ** If we want to only impose $p>2$, then we either need
%to get more moments for $\Delta$, or just be satisfied with 4 moments for
%$\Delta$, and improve on the first term in \eqref{PBn} instead. Note that
%the second term there is OK with $p=2$. **}

\begin{proof}[Proof of Lemma \ref{Delta-moments}]
The measurability statement is immediate from the definition of
$\Delta$. Concerning the moments,
first note that the moment hypothesis (M) implies that $g$ is
uniformly bounded by $2M$. Second, note that combining the non-nestling
hypothesis (N) with the moment hypothesis (M) and H\"older's inequality,
one has
\[\delta\leq\sum_{x\geq1}x\pi_{0x}\leq
(1-\pi_{00})^{(p-1)/p}M,\]
which says that $\pi_{00}$ is uniformly bounded away from $1$.
Third, note that if we denote by $P_{i,j}^\w$ the
process of two independent walkers in the same environment $\w$, one starting at $i$
and the other at $j$, then
\[a_i(\w)=(1-\pi_{ii}(\w))^{-1}E_{0,1}^\w(\one_{V_i}-\one_{\bar V_i}),\]
where
$\bar V_i$ is the event corresponding to the walker starting at $1$. If we denote
by $L$ the first common point of the paths $\{X_n: n\geq 0\}$ and
$\{\bar{X}_n: n\geq 0\}$, then we have
\begin{align*}
|a_i(\w)|&\leq(M/\delta)^{p/(p-1)}\Big|\sum_{j=1}^i P_{0,1}^\w(L=j)
E_{j,j}^\w(\one_{V_i}-\one_{\bar V_i})\Big|\\
&\qquad+(M/\delta)^{p/(p-1)}
\Big|E_{0,1}^\w(\one_{V_i}-\one_{\bar V_i},L>i)\Big|\\
&\leq (M/\delta)^{p/(p-1)}P_{0,1}^\w(L>i).
\end{align*}
The coefficients $P_0^\w(V_i)$ and $P_1^\w(V_i)$ are functions of
$\w_0,\dotsc,\w_{i-1}$, hence independent of $(1-\pi_{ii}(\w))^{-1}g(T_i\w)$.
By \eqref{Pinftyformula} and the definition of $v$ we have
\[\E[g/(1-\pi_{00})]=\E_\infty(g)E_0(\si_1)=0.\]
 Therefore, under the measure $\P$, 
$\Delta_n=\sum_{i=0}^n a_i(\w) g(T_i\w)$ is a martingale relative to the
filtration $\si(\w_0,\dotsc,\w_n)$. Then the Burkholder-Davis-Gundy inequality
%\cite[Th. 48 on p. 193]{protter}
\cite[Theorem 3.2]{burk} gives us
\begin{align*}
\sup_n \E[|\Delta_n|^m]&\leq C\,\sup_n
\E\left[\Big|\sum_{i=0}^n |a_i|^2\Big|^{m/2}\right]
\leq C^{'}\, \E\left[\Big|\sum_{i\geq0} P_{0,1}^\w(L>i)\Big|^{m/2}\right]\\
&\leq C^{'}\, \E[|E_{0,1}^\w(L)|^{m/2}] \leq C^{'}\, E_{0,1}[L^{m/2}].
\end{align*}
In the last stage above we can take  
$P_{0,1}=P_0\otimes P_1$ to be  the process of two
 walkers with independent annealed distributions. 
This works because 
 the two walkers do not meet until they visit site  $L$ for the first time. 
We can apply Lemma \ref{renewal-lm} from the appendix to $L$, by taking 
$Y_i=X_{\si_i}-X_{\si_{i-1}}$ with distribution 
\[
P[Y_i=k] = \E[(1-\pi_{00})^{-1}\pi_{0k}]\,, \quad k\geq 1.
\]  
Bound \eqref{5} from Lemma \ref{exponential} gives 
$E(Y_i^{\bar p})<\infty$ for any $\bar p<p$. 
Then  Lemma \ref{renewal-lm}  tells us  that $E_{0,1}(L^{\bar p-1})<\infty$
and  we can take $m=2\bar p-2$ in the bound above. 
%Now, we can replace the two annealed walks by ones with transition vector
%$\{(1-\\pi_{00})^{-1}\E\pi_{0z}\}_{z>0}$ without changing the
%distribution of $L$.
%Lemma \ref{renewal-lm} then tells us that $E_{0,1}(L^{2[p/2]})<\infty$.

The integrability of $L$ and the earlier bound on $|a_i(\w)|$ imply that the 
series \eqref{Delta} defining $\Delta$ converges absolutely,
$\P$-almost surely. 
By the moment bound above and the martingale convergence theorem, 
$\Delta$ is also the $L^m(\P)$-limit of $\Delta_n$. This gives 
$\E\lvert\Delta\rvert^m<\infty$ and 
 $\E[\Delta]=0$.
\end{proof}

\begin{remark}
A place where it should be possible
to improve the moment hypotheses needed for the argument is the moment bound
on $\Delta$. For example, one has
\[
\sum_{i\geq 0} \abs{a_i}^2 \leq \sum_{i\geq0} P_{0,1}^\w(L>i)^2
%= \sum_i P^\w( L \land L' > i )   [independent copy L' of L]
=E^\w( L\land \bar L),
\]
where $\bar L$ is an independent copy of $L$.
Then for the $m$-th moment  of $\Delta_n$  we get the bound
\[
\E[ E^\w( L\land \bar L) ^{m/2} ]
\leq \E[ E^\w( L^{m/4} \land {\bar L}^{m/4}) ^{2} ]
\leq \E[ E^\w( L^{m/4} )^2 ].
\]
Now $m$ is divided by 4 instead of 2, but the drawback is that
one needs to bound the quenched moment of $L$.
\end{remark}

Define $\chi(0,\w)=0$, and for $x\geq1$ \[\chi(x,\w)=\sum_{y=0}^{x-1}\Delta(T_y\w).\]
One can check that $a_0=1/(1-\pi_{00})$ and, for $i\geq1$,
\begin{align*}
\sum_{j=0}^{i} a_j(T_{i-j}\,\w)\sum_{y>i-j}\pi_{0y}(\w)=0.
%\label{a-aux}
\end{align*}
%let V_{ij}=P_i^\w(V_j). then
%a_j(T_{i-j}\w)=(1-\pi_{ii})^{-1}V_{i-j,i}-V_{i+1-j,i}.
%change indices k=i-j.  call the cumulative distribution function F_k,
%so that F_{k-1}-F_k=\pi_{0k}).  now, use summation by parts.
With the help of the above formula one can check that 
\begin{align}
\label{Echi=g}
E_0^\w[\chi(X_1,\w)]=g(\w).
\end{align}
%Indeed
%\begin{align*}
%\Detla(\w)&=g(\w)-\sum_{i\geq1}\sum_{j=0}^{i-1}a_j(T_{i-j}\w)g(T_i\w)\sum_{y>i-j}\pi_{0y}(\w)\\
%&=g(\w)-\sum_{j\geq0}\sum_{k\geq1}a_j(T_k\w)g(T_{k+j}\w)\sum_{y>k}\pi_{0y}(\w)\\
%&=g(\w)-\sum_{k\geq1}\Delta(T_k\w)\sum_{y>k}\pi_{0y}(\w).
%\end{align*}
%Exchanging the order of summation gives then
%\[\sum_{y>0}\pi_{0y}(\w)\chi(y,\w)=g(\w).\]
%But now,
%\begin{align*}
%E_{X_n}^\w[\chi(X_{n+1},\w)]&=E_{X_n}^\w\left[\sum_{y=0}^{X_{n+1}-1}\Delta(T_y\w)\right]\\
%&=\chi(X_n,\w)+E_{X_n}^\w\left(\sum_{y=X_n}^{X_{n+1}-1}\Delta(T_y\w)\right)\\
%&=\chi(X_n,\w)+E_0^{\w}[\chi(X_1,\w)]\Big|_{\w\leftrightarrow T_{X_n}\w}\\
%&=\chi(X_n,\w)+g(T^{X_n}\w).
%\end{align*}
%
This implies that, for $\P$-a.e. $\w$,
\begin{align}
\label{Mn}
M_n=X_n-nv-\chi(X_n,\w)
\end{align}
is a $P_0^\w$-martingale.
Using the invariance and ergodicity
 of $\P_\infty$ under the process $(T_{X_n}\w)$,
and the ergodic theorem, one can check that $M_n$ satisfies the conditions of the
martingale invariance principle; see, for example, Theorem (7.4) of Chapter 7
of \cite{durrett}.
% Basically, $E_0|\chi(X_1,\w)|^2\leq E_0 X_1^2\E|\Delta|^2.
So the distribution of $(n^{-1/2}M_{[nt]})_{t\geq0}$
under $P_0^\w$ converges weakly to a Brownian motion with
diffusion coefficient $\E_\infty E_0^\w(M_1^2)$, for $\P_\infty$-a.e.\ $\w$.
Since everything in this statement is $\kS_0$-measurable, the same invariance principle
holds for $\P$-a.e.\ $\w$.

To handle the random summation limit  in $\chi(X_n,\w)$
we decompose 
\begin{align}
\label{chi}
\chi(X_n,\w)=Z_n(\w)+R_n
\end{align}
where $Z_n(\w)=\sum_{y=0}^{[nv]}\Delta(T_y\w)$ and $R_n$ is the 
remaining part.
%We have
%\begin{align}
%\label{decomp}
%X_n-E_0^\w[X_n]=M_n+R_n-E_0^\w[R_n].
%\end{align}

We can now outline the strategy for proving the invariance 
principles. From 
\eqref{Mn} and \eqref{chi} we write 
\be
E_0^\w(X_n)-nv=Z_n(\w)+E_0^\w(R_n)\label{EwXn} 
\ee
and then 
\be 
X_n-E_0^\w(X_n)=M_n+R_n-E_0^\w(R_n).
\label{XnMn}
\ee
With \eqref{XnMn} 
the quenched invariance principle  for the processes 
$\Btil_n(t)=n^{-1/2}\{X_{[nt]}-E_0^\w(X_{[nt]})\}$ 
 follows from the invariance
principle of $M_n$  observed above, after we show that the 
errors $R_n$ and $E_0^\w(R_n)$ are negligible.
For the quenched mean process $n^{-1/2}\{E^\w_0(X_{[nt]})-ntv\}$
we show that the process $Z_n$ 
is close to a martingale under $\P$, and apply the martingale
invariance principle again.  

As a preparatory step we state an annealed bound on deviations
of the walk.

\begin{lemma} For every  $\bar p\in[2,p)$ 
there exists a constant $C=C(\bar p)<\infty$ such that, 
for $h>0$ and $n\geq 1$, 
\be
P_0(|X_n-nv|> h)\leq Ch^{-\bar{p}} n^{\bar{p}/2}.
\label{ann-mdp}
\ee
\label{ann-mdp-lm}
\end{lemma}

\begin{proof}[Proof of Lemma \ref{ann-mdp-lm}]
Recall the definition $K(n)=\max\{j: \si_j\leq n\}$. 
In one dimension 
$\si_{K(n)+1}=\inf\{ k>n: X_k>X_n\}$ because $X_n$ does not jump 
during any time interval $\si_i\leq n<\si_{i+1}$.
% which makes
%$K(n)+1$ into a stopping time for the process $X_{\si_k}-\si_kv$. 
%For a walk
%in one dimension, $X_{\si_{K(n)}}=X_n$ because $X_n$ does not jump 
%during any time interval $\si_i\leq n<\si_{i+1}$.
\begin{align*}
&P_0(\,|X_n-nv|>h)\leq 
P_0(\,|X_n-X_{\si_{K(n)+1}}|>h/3) \\
&\qquad + P_0(\,|X_{\si_{K(n)+1}}-\si_{K(n)+1}v|>h/3) +
P_0(\si_{K(n)+1}- n > h |v|^{-1}/3)\\
&\leq E_0 P^\w_{X_n}(X_{\si_1}-X_0 >h/3)
+ P_0\bigl(\, \max_{k\leq n+1} |X_{\si_k}-\si_kv|>h/3\bigr)\\
&\qquad +
E_0 P^\w_{X_n} (\si_1 >h |v|^{-1}/3).
\end{align*}
For the last inequality we restarted the walk at time $n$ 
for the first and last probability, and used $K_n\leq n$ for the 
middle probability. 
>From \eqref{5} we get a bound on the first probability:
\[
E_0 P^\w_{X_n}(X_{\si_1}-X_0 >h/3)
\leq Ch^{-\bar{p}}\leq Ch^{-\bar{p}} n^{\bar{p}/2}.  
\]
The third probability is handled similarly via \eqref{4}. 
For the middle probability note that $X_{\si_k}-\si_kv$ is a sum
of mean zero i.i.d.~random variables with finite $\bar p$-th moment. 
 Use first Doob's inequality,  then the Burkholder-Davis-Gundy inequality
\cite[Theorem 3.2]{burk}, and finally the Schwarz inequality: 
\begin{align*}
&P_0\bigl(\, \max_{k\leq n+1} |X_{\si_k}-\si_kv|>h/3\bigr)
\leq Ch^{-\bar{p}} E_0\bigl[ \,|X_{\si_{n+1}}-\si_{n+1}v|^{\bar{p}}\;\bigr]\\
&\leq Ch^{-\bar{p}} E_0\biggl[ \,\biggl\lvert \sum_{k=1}^{n+1} 
(X_{\si_{k}}-X_{\si_{k-1}}-\si_{k}v+\si_{k-1}v)^2\biggr\rvert^{\bar{p}/2}\;\biggr]\\
&\leq Ch^{-\bar{p}} (n+1)^{\bar{p}/2} E_0\biggl[ \, \frac1{n+1}\sum_{k=1}^{n+1} 
\lvert X_{\si_{k}}-X_{\si_{k-1}}-\si_{k}v+\si_{k-1}v \rvert^{\bar{p}}\;\biggr]
 \leq C^{'}h^{-\bar{p}} n^{\bar{p}/2}.
\end{align*}
Inequality \eqref{ann-mdp} has been verified. 
\end{proof}

We turn to the first main task, bounding the errors 
$R_n$ and $E_0^\w(R_n)$.
Fix $\e>0$ and  $\eta\in(0,1/2-1/p)$, and let $\delta_n=n^{-\eta}$. 
The strategy is to  develop summable deviation estimates 
for the errors and apply Borel-Cantelli. Summability of the
estimates  will 
require that $\eta$ satisfies several conditions simultaneously.
Then we show that for %$p>\tfrac12(5+\sqrt{17})$  
$p>p_0$
a single choice of $\eta$ can satisfy
all the requirements. 
%Observe that under $P_0$ the walk is just a homogeneous random walk with
%independent increments.
Write
\begin{align}
P_0(|R_n|>\e\sqrt{n})
%\leq \P(|X_n-nv|>n\delta_n) +
%\P(|X_n-nv|\leq n\delta_n,|E^\w_0[R_n]|>\e\sqrt{n})
\leq P_0(|X_n-nv|>n\delta_n) +
\P(\cB_n^\e)+\P(\cC_n^\e)
\label{inequa}
\end{align}
where
\[\cB_n^\e=\left\{\max_{0\leq i\leq n\delta_n}
\abs{\sum_{j=0}^i\frac{\Delta(T_{-j}\w)}{\sqrt n}}>\frac\e2\right\}
\ \text{ and }\
\cC_n^\e=\left\{\max_{0\leq i\leq n\delta_n}
\abs{\sum_{j=i}^{[n\delta_n]}\frac{\Delta(T_{-j}\w)}{\sqrt n}}>\frac\e2\right\}.\]

If we pick $\bar{p}<p$ close enough to $p$ so that 
\be
p^{-1}<\bar{p}^{\,-1}<\tfrac12-\eta
\label{eta-cond-1}
\ee
 to ensure $\bar{p}(\tfrac12-\eta)>1$,
Lemma \ref{ann-mdp-lm} shows that the first term on the right-hand side of
\eqref{inequa} is summable:
\be
\sum_n P_0(\,|X_n-nv|>n\delta_n) \leq C \sum_n n^{-\bar{p}(\tfrac12-\eta)} < \infty.
\label{ann-mdp-2}
\ee
%\begin{align*}
%P_0(|X_n-nv|>n\delta_n)=\Ord(\frac{n^{p/2}}{(n\delta_n)^p}).
%\end{align*}

To control the other two terms in \eqref{inequa}, define
\begin{align*}
H_n=\sum_{k\geq0}\Delta(T_{-n+k}\w)\sum_{y>k}\E_\infty\pi_{0y}.
\end{align*}
For the next calculations,
recall that 
\[\frac{d\P_\infty}{d\P}=\frac{(1-\pi_{00})^{-1}}{\E[(1-\pi_{00})^{-1}]}
\quad\text{and}\quad
v=\E_\infty D =\E_\infty \sum_{k\geq 0}\sum_{y>k}\pi_{0y}. 
\]
Then, for $m=2\bar p-2$,
\begin{align*}
\E[|H_n|^m]&=v^m\E\left[\Big|v^{-1}\sum_{k\geq0}\Delta(T_k\w)
\sum_{y>k}\E_\infty\pi_{0y}\Big|^m\right]\\
&\leq v^{m-1}\sum_{k\geq0}\E[|\Delta|^m]\sum_{y>k}\E_\infty\pi_{0y}\\
&=v^m\E[|\Delta|^m]<\infty.
\end{align*}
Also, under $\P$, $H_n$ is a martingale difference relative to the 
filtration $\kS_{-n}$. To see that, it is enough to show that 
$\E(H_0|\kS_1)=0$, since $H_n=T_{-n}H_0$. And indeed:
\begin{align*}
\E[(1-\pi_{00})^{-1}]\E[H_0|\kS_1]&=\E[\Delta|\kS_1]+
\sum_{k\geq1}\Delta(T_k\w)\sum_{y>k}\E[(1-\pi_{00})^{-1}\pi_{0y}]\\
&=\E\Big[(1-\pi_{00})^{-1}\sum_{k\geq0}\Delta(T_k\w)\sum_{y>k}\pi_{0y}\Big|\kS_1\Big]\\
&=\E[(1-\pi_{00})^{-1}g]=0,
\end{align*}
where we have used \eqref{Echi=g} and $\E_\infty g=0$.
Moreover,
\begin{align*}
\sum_{k=1}^n (v\Delta(T_{-k}\w)-H_k)
&=\sum_{j\geq0}\sum_{y>j}\E_\infty\pi_{0y}\sum_{k=1}^n(\Delta(T_{-k}\w)-\Delta(T_{-k+j}\,\w))\\
&=\sum_{y>1}\E_\infty\pi_{0y}\sum_{j=1}^{y-1}\sum_{k=1}^n(\Delta(T_{-k}\w)-\Delta(T_{-k+j}\,\w))\\
&=\sum_{y>1}\E_\infty\pi_{0y}\sum_{j=1}^{y-1}\sum_{k=0}^{j-1}(\Delta(T_{k-n}\w)-\Delta(T_k\w)),
\end{align*}
where the last line was obtained by adding (resp.\ subtracting)
the necessary terms in the third summation when $j\geq n$ (resp.\ $j<n$).
Thus
\begin{align}
\label{Hn}
v\sum_{k=0}^n\Delta(T_{-k}\w)=\sum_{k=1}^n H_k+ U_n,
\end{align}
with
\begin{align*}
U_n(\w)&=v\Delta(\w)+A(T_{-n}\w)-A(\w),\\
A(\w)&=\sum_{y>1}\E_\infty(\pi_{0y})\sum_{j=0}^{y-2}(y-j-1)\Delta(T_j\w).
\end{align*}
Now, similarly to the moment bound above for $H_n$, 
\begin{align*}
\E[|A|^m]\leq\E\left[\Big|\sum_{j\geq0}|\Delta(T_j\w)|\sum_{y>j}y\E_\infty\pi_{0y}\Big|^m\right]
\leq [\E_\infty E_0^\w X_1^2]^m\E[|\Delta|^m].
\end{align*}
Then, using Doob's inequality, followed by the Burkholder-Davis-Gundy
inequality \cite[Theorem 3.2]{burk}, 
and then by Jensen's inequality on the first
term below, and Chebyshev's inequality on the second term, one has
%Given my discrete time mg M_n with increments
%H_i=M_i-M_{i-1}, define the cadlag mg  M_t=M_{[t]}.
%the quadratic variation is
%[M]_n = \sum_{i=1}^n H_i^2
%so BDG gives
%E[ \max_{k\leq n} |M_k|^p ] \leq  C E[ (\sum_{i=1}^n H_i^2 )^{p/2}]
%\leq  c n^{p/2} E[ n^{-1}\sum_{i=1}^n H_i^p  ]
%And now bounded p moment for H_i
\begin{align}
\P(\cB_n^\e)&\leq \P\left(\max_{i\leq n\delta_n}\abs{\sum_{k=1}^i H_k}>
\frac{v\e\sqrt n}4\right)+
n\delta_n\max_{i\leq n\delta_n}\P\left(|U_i|>\frac{v\e\sqrt n}4\right)
\nonumber\nn\\
&\leq Cn^{-m/2} \E \biggl\lvert \sum_{k=1}^{[n\delta_n]} H_k\biggr\rvert^{m}
+ C\delta_nn^{1-m/2}\nn\\
&\leq Cn^{-m/2} \E \biggl\lvert \sum_{k=1}^{[n\delta_n]} H_k^2\biggr\rvert^{m/2}
+ C\delta_nn^{1-m/2}\nn\\
&\leq Cn^{-m/2}  (n\delta_n)^{m/2-1} \E  \sum_{k=1}^{[n\delta_n]} \lvert H_k\rvert^m 
+ C\delta_nn^{1-m/2}\nn\\
&\leq C(\delta_n^{m/2}+\delta_n n^{1-m/2}). \label{PBn}
\end{align}
%The application of Jensen's inequality with $m/2$ requires $m\geq 2$,
%which in turn needs $\bar p\geq 2$. 

%{\bf ** The second term in \eqref{PBn} is OK once $m\geq4$. This is OK, even with
%$p=2$. For the first term and the term controlling $P_0(|X_n-nv|>n\delta_n)$, we
%need $m>4p/(p-2)$. If the highest we can take is $m=4[p/2]$, then the above
%works once $p>4$. This is weaker than the restriction
%needed to workout the $E_0^\w R_n$ part below.
%
%Note also that with $p$ close to $2$, we need a very large $m$. This means that we should
%either improve our moment bounds on $\Delta$, or get $m=2p$ instead of only
%integer $m$, and improve on either the first bound in \eqref{PBn}
%or the LDP bound on $P_0(|X_n-nv|>n\delta_n)$ (to allow just $\eta<1/2$). **}

%If $p\geq 5$, we can choose $\eta\in(0,1/2-1/p)$ so that 
%$\sum_n\P(\cB_n^\e)<\infty$.
The same argument works  also for $\P(\cC_n^\e)$,
by first writing 
\[
\biggl\lvert \sum_{j=i}^{[n\delta_n]}\Delta(T_{-j}\w)\biggr\rvert
\leq
\biggl\lvert \sum_{j=0}^{[n\delta_n]}\Delta(T_{-j}\w)\biggr\rvert
+
\biggl\lvert \sum_{j=0}^{i-1}\Delta(T_{-j}\w)\biggr\rvert.
\]
In order to have the Borel-Cantelli Lemma imply that
$R_n/\sqrt n$ goes to $0$ $\P$-a.s.  and,
therefore, also 
\begin{align}
\label{Rn-small}
\max_{k\leq n}|R_k|/\sqrt n\hbox{ goes to 0, }\P\hbox{-a.s.}
\end{align}
we need 
\be
\sum_n (\delta_n^{m/2}+\delta_n n^{1-m/2})
=\sum_n (n^{-\eta m/2} + n^{-\eta+1-m/2}) <\infty.
\label{eta-cond-2}
\ee

Next, write
\begin{align*}
\P(|E_0^\w[R_n]|>\e\sqrt{n})\leq
&\P(|E_0^\w[R_n,X_n-nv>n\delta_n]|>\e\sqrt n/4)\\
&\quad+\P(|E_0^\w[R_n,nv-X_n>n\delta_n]|>\e\sqrt n/4)\\
&\quad+\P(\cB_n^{\e/2})+\P(\cC_n^{\e/2})
%+\P(T^{[nv]}\w\in\cB_n)+\P(T^{[nv]+[n\delta_n]}\w\in\cC_n).
%\label{inequa}
\end{align*}
The last two terms can be controlled as in \eqref{PBn}.
The first two are similar to each other,
so we work only with  the first one. To this end, we have
\begin{align*}
&\P(|E_0^\w[R_n,X_n-nv>n\delta_n]|>\e\sqrt n/4)\\
&\qquad\leq\frac4{\e\sqrt n}\sum_{x>nv+n\delta_n}
\E\left[P_0^\w(X_n=x)\sum_{y=[nv]}^{x-1}|\Delta(T_y\w)|\right]\\
&\qquad\leq\frac4{\e\sqrt n}\sum_{x>nv+n\delta_n}
\E[P_0^\w(X_n=x)^{m/(m-1)}]^{1-1/m}\,\E\left[\Big|\sum_{y=[nv]}^{x-1}|\Delta(T_y\w)|\Big|^m\right]^{1/m}\\
&\qquad\leq\frac4{\e\sqrt n}\sum_{x>nv+n\delta_n}
P_0(X_n=x)^{1-1/m}\,\E\left[(x-[nv])^{m-1}\sum_{y=[nv]}^{x-1}|\Delta(T_y\w)|^m\right]^{1/m}\\
&\qquad\leq\frac{4\E[|\Delta|^m]^{1/m}}{\e\sqrt n}
\sum_{x>n\delta_n}x\,P_0(X_n-[nv]=x)^{1-1/m}\\
&\qquad \leq Cn^{-1/2+\bar{p}(1-1/m)/2}\sum_{x>n\delta_n}x^{1-\bar{p}(1-1/m)}
\leq 
%Cn^{-1/2+\bar{p}(1-1/m)/2+2-\bar{p}(1-1/m)-\eta(2-\bar{p}(1-1/m))} =
Cn^{3/2-2\eta -\bar{p}(1/2-\eta)(1-1/m)},
%&\qquad\leq\frac{4\E[|\Delta|^m]^{1/m}}{\e\sqrt n}
%\sum_{x>n\delta_n}\frac{\Ord(n^{(1-1/m)p/2})}{x^{(1-1/m)p-1}}=\Ord(n^{-b}),
\end{align*}
where the next to last inequality came from Lemma \ref{ann-mdp-lm}. 
%with $2mb=2(p(m-1)-2m)(1-\eta)-p(m-1)+m$. If $p\geq7$,
%then one can choose $\eta\in(0,1/2-1/p)$ so that  $b>1$.
%
%{\bf ** If we improve on this bound, without doing anything else, then we can
%have $p>4$ as a restriction coming from \eqref{PBn}. **}
%
Again, the Borel-Cantelli Lemma implies that
\begin{align}
\label{ERn-small}
\max_{k\leq n}|E_0^\w(R_k)|/\sqrt n\hbox{ goes to 0, }\P\hbox{-a.s.}
\end{align}
if 
\be
3/2-2\eta -\bar{p}(1/2-\eta)(1-1/m)<-1.
\label{eta-cond-3}
\ee

Recall the definition \eqref{def-p0} of $p_0$. 
Now we observe that 
$p>p_0$
%$p>\tfrac12(5+\sqrt{17})\approx 4.6$ 
enables us to make all the Borel-Cantelli
arguments simultaneously valid.  As throughout this section,
$m=2\bar p-2$. 

\begin{lemma} Suppose %$p>\tfrac12(5+\sqrt{17})$. 
$p>p_0$.
Then we can choose 
$\eta\in(0, 1/2-1/p)$ so that there exists  $\e_1>0$ 
such that, for all $\bar p\in(p-\e_1,p)$, requirements
{\rm\eqref{eta-cond-1}}, {\rm\eqref{eta-cond-2}},  and
{\rm\eqref{eta-cond-3}} all hold,   and thereby the 
conclusions {\rm\eqref{Rn-small}} and {\rm\eqref{ERn-small}}
hold. 
\label{eta-lemma}
\end{lemma}

\begin{proof}[Proof of Lemma \ref{eta-lemma}]
The first requirement is easy: for any $\eta\in(0,1/2-1/p)$
there is a whole range $(p-\e_2,p)$ of $\bar p$-values that satisfy
\eqref{eta-cond-1}, and forcing $\bar p$ closer to $p$ cannot
violate this condition. 

The second condition \eqref{eta-cond-2} has two parts. 
The first one requires $1<\eta m/2=\eta(\bar p-1)$, or
 $\eta>(\bar p-1)^{-1}$. 
%This is where the condition on $p$ comes from:  
Observe that $\tfrac12(5+\sqrt{17})\approx4.6$  is the larger root of
$(p-1)^{-1}=1/2-1/p$. Hence for $p>\tfrac12(5+\sqrt{17})$  
 and $\bar p$  close enough to $p$, there is room to choose 
$\eta$ so that $(\bar p-1)^{-1}<\eta<1/2-1/p$. 

The second part of \eqref{eta-cond-2} requires 
$-\eta+1-m/2<-1$ which is satisfied for any $\eta>0$ if $\bar p\geq 3$.
With $p> \tfrac12(5+\sqrt{17})$ and $\bar p$ close enough to $p$ 
this is guaranteed. 

Finally we consider requirement \eqref{eta-cond-3}. 
Take first $\eta=(\bar p-1)^{-1}$ which was identified as a lower
bound for $\eta$ two paragraphs above.
By getting rid of denominators, observe  that \eqref{eta-cond-3}  
is satisfied iff
\[
 (2\bar p-2)(5\bar p-9)-\bar p(\bar p-3)(2\bar p-3) <0.
\]
By calculus we see that over the interval $(4,\infty)$, this happens iff
$p$ is larger than $p_0=\frac{19}6+\frac{\sqrt{139}}3\cos(\frac13\arccos(\frac{1504}{139^{3/2}}))$.
%which by calculus is true for example for $\bar p>4$.
Since the inequality above is strict, it remains in force if $\eta>(\bar p-1)^{-1}$
is close enough to $(\bar p-1)^{-1}$. 

To summarize, as long as %$p>\tfrac12(5+\sqrt{17})$, 
$p>p_0$,
a pair $(\bar p,\eta)$ can be 
found to satisfy all the requirements. 
\end{proof} 

We are ready to prove  Theorem \ref{qclt-1d}. 
Part (a)  follows 
from \eqref{XnMn},
and the martingale invariance principle applied to $M_n$, 
because \eqref{Rn-small} and \eqref{ERn-small} show
 that $R_n-E_0^\w(R_n)$ is a negligible error.
Formula \eqref{kappa_q} says simply that the diffusion coefficient is
$\E_\infty E_0^\w |M_1|^2$.

For part (b), combine \eqref{Hn} with \eqref{EwXn} to write
\be
\begin{split}
E_0^\w(X_n)-nv&=Z_n+E_0^\w(R_n)\\
&=v^{-1}\sum_{k=1}^{[nv]}H_k(T_{[nv]}\w)+
v^{-1}U_{[nv]}(T_{[nv]}\w)+E_0^\w(R_n).
\end{split} \label{q-decomp}
\ee
By the calculation in \eqref{PBn}, one can neglect $U_n$. By \eqref{ERn-small},
one can also neglect
$E_0^\w(R_n)$. Then, using the ergodicity of $\P$, one can check that the martingale
$\sum_{k=1}^nH_k$ satisfies the conditions
of the martingale invariance principle 
\cite[Section 7.7]{durrett}. %since H_n=T_nH_0
This implies (b), short of identifying the diffusion coefficient.

If we define $\kappa_m^2$ as the limiting diffusion coefficient 
for the process 
$n^{-1/2}\{E_0^\w(X_{[nt]})-ntv\}$, then \eqref{q-decomp} implies
 $\E(H_0^2)=\kappa_m^2 v$.
Then, to calculate $\E(H_0^2)$, use 
 the fact that the partial sums of the infinite sum
below form a martingale relative
to the filtration $\sigma(\w_0,\cdots,\w_n)$: 
\begin{align*}
\E[H_0^2]&=\E\left[\Big|
\sum_{i\geq0}g(T_i\w)\sum_{k=0}^i a_{i-k}(T_k\w)\sum_{y>k}\E_\infty\pi_{0y}
\Big|^2\right]\\
&=\E[g^2/(1-\pi_{00})^2]
%\sum_{i\geq0}\E\left[\Big|
%\sum_{k=0}^i a_{i-k}(T_k\w)\sum_{y>k}\E_\infty\pi_{0y}
%\Big|^2\right]\\
%=\E[|E_0^\w(X_{\si_1})-v E_0^\w(\si_1)|^2]
\sum_{i\geq0}\E\left[\Big|
\sum_{k=0}^i (P_k^\w(V_i)-P_{k+1}^\w(V_i))\sum_{y>k}\E_\infty\pi_{0y}
\Big|^2\right]\\
&=\E[|E_0^\w(X_{\si_1})-v E_0^\w(\si_1)|^2]\sum_{i\geq0}\E\left[\Big|
P_0^\w(V_i)-\sum_{k=0}^i P_k^\w(V_i)\E_\infty\pi_{0k}
\Big|^2\right].
\end{align*}
A comparison with \eqref{kappa_m} shows that
this completes the proof of part (b). 

As for part (c), use \eqref{q-decomp} to write
\[\frac{X_n-nv}{\sqrt n}=\frac{X_n-E_0^\w(X_n)}{\sqrt n}+
\frac{v^{-1}U_{[nv]}(T_{[nv]}\w)+E_0^\w(R_n)}{\sqrt n}+
\frac{\sum_{k=1}^{[nv]}H_k(T_{[nv]}\w)}{v\sqrt n}.\]
The first term above is tight under $P_0^\w$, for $\P$-a.e.\ $\w$,
by part (a). We have shown that the second term
goes to $0$ in $\P$-probability. For the last term, we have also shown that
it converges to a Gaussian under $\P$.
This implies that
\[\P\left(\lsup_{n\to\infty}\frac{\sum_{k=1}^n H_k(T_n\w)}{\sqrt n}\geq y\right)>0,\quad\forall y.\]
Since the above is a tail event and $\P$ is product, the third term cannot be bounded,
for a fixed $\w$, unless $0=\E(H_0^2)=v\kappa_m^2$.
%$0=\E(H_0^2)=\E(g^2/(1-\pi_{00})^2)(\E D)^2/E_0(X_{\si_1})$.
But, according to part (a) of Proposition \ref{km+kq},
this is not allowed by the non-degeneracy condition.
\end{proof}

We conclude this section by proving 
Proposition {\rm\ref{km+kq}} about the diffusion coefficients. 

\begin{proof}[Proof of Proposition {\rm\ref{km+kq}}]
To prove part (a) one notices that in the sum in \eqref{kappa_m} 
 the $i=0$ term is positive unless the walk never
leaves $0$.  But as observed in the proof of  
 Lemma \ref{Delta-moments}, $\pi_{00}$ is bounded away from $1$. 
Therefore, $\kappa_m=0$ if, and only if, $(1-\pi_{00})^{-1}g=0$, $\P$-a.s.
Part (a) follows.

To prove part (b) we will distinguish two cases. First, if $\P(\pi_{00}>0)>0$, then
\[%\E_\infty E_0^\w\Big|X_1-v-\sum_{z=0}^{X_1-1}\Delta(T_z\w)\Big|^2
\kappa_q^2\geq v^2(E_0\si_1)^{-1}\E\Big[\frac{\pi_{00}}{1-\pi_{00}}\Big]>0.\]

On the other hand, if $\pi_{00}=0$, $\P$-a.s., and if
there exist two distinct points $x$ and $y$ such that $\E\pi_{0x}\pi_{0y}>0$, then
\[\E\Big[\prod_{i=0}^{y-1}\pi_{ix,(i+1)x}\prod_{i=0}^{x-1}\pi_{iy,(i+1)y}\Big]>0\]
and if $X_1-v-\sum_{z=0}^{X_1-1}\Delta(T_z\w)=0$, $P_0$-a.s., then
\[\sum_{z=0}^{xy-1}(1-\Delta(T_z\w))\]
is simultaneously equal to
\[\sum_{i=0}^{y-1}\sum_{z=0}^{x-1}(1-\Delta(T_{z+ix}\w))=yv\]
and
\[\sum_{i=0}^{x-1}\sum_{z=0}^{y-1}(1-\Delta(T_{z+iy}\w))=xv\]
which implies that $v=0$ and is a contradiction.

The proof of part (c) uses  uniform integrability that
follows from having $p>2$. 
%%%%%%%%The details are left to the reader.
Start with 
\be
E_0\biggl[\,\biggl(\frac{X_n-nv}{\sqrt{n}}\biggr)^2\,\biggr]
=
E_0\biggl[\,\biggl(\frac{X_n-E_0^\w(X_n)}{\sqrt{n}}\biggr)^2\,\biggr]
+
\E\biggl[\,\biggl(\frac{E^\w_0(X_n)-nv}{\sqrt{n}}\biggr)^2\,\biggr].
\label{var-decomp-2}
\ee
Let $2<\tilde p<\bar p<p$. Estimate \eqref{ann-mdp} implies
that 
\[
\sup_n E_0\biggl[\,\biggl\lvert
\frac{X_n-nv}{\sqrt{n}}\biggr\rvert^{\tilde p}\,\biggr]
<\infty,
\]
which together with the annealed CLT (Theorem \ref{a-clt}) implies
\[
\lim_{n\to\infty} 
E_0\biggl[\,\biggl(\frac{X_n-nv}{\sqrt{n}}\biggr)^2\,\biggr]
=\frac{E_0(\,\lvert X_{\si_1}-\si_1v\rvert^2\,)}{E_0(\si_1)}.
\]
Similarly on the right-hand side of \eqref{var-decomp-2} we get
convergence to $\kappa_q^2+\kappa_m^2$ by observing that 
uniform integrability follows from that already proved:
\[
 \E\biggl[\,\biggl\lvert
\frac{E^\w_0(X_n)-nv}{\sqrt{n}}\biggr\rvert^{\tilde p}\,\biggr]
=
 \E\biggl[\,\biggl\lvert
E^\w_0\biggl(\frac{X_n-nv}{\sqrt{n}}\biggr) \biggr\rvert^{\tilde p}\,\biggr]
\leq
 E_0\biggl[\,\biggl\lvert
\frac{X_n-nv}{\sqrt{n}}\biggr\rvert^{\tilde p}\,\biggr]
\]
and then 
\[
E_0\biggl[\,\biggl\lvert
\frac{X_n-E^\w_0(X_n)}{\sqrt{n}}\biggr\rvert^{\tilde p}\,\biggr]
\leq
2^{\tilde p} E_0\biggl[\,\biggl\lvert
\frac{X_n-nv}{\sqrt{n}} \biggr\rvert^{\tilde p}\,\biggr]
+  2^{\tilde p}
 \E\biggl[\,\biggl\lvert
\frac{E^\w_0(X_n)-nv}{\sqrt{n}}\biggr\rvert^{\tilde p}\,\biggr].
\]

To prove part (d), recall that $K(n)=\max\{j:\sigma_j\leq n\}$
and observe that
\begin{align}
\E[|E_0^\w(X_n)-nv|^2]&=\E\biggl[\bigg|\sum_{k=0}^{n-1}
E_0^\w[g(T_{X_k}\w)]\bigg|^2\biggr]\nn\\
&=\E\biggl[\bigg|E^\w_0\sum_{k=0}^{\si_{K(n)+1}-1}  g(T_{X_k}\w)\bigg|^2\biggr] +\ord(n). \label{error-1}
\end{align}
The bound 
  $\ord(n)$ on the error above comes because of ergodicity of $\P_\infty$, 
boundedness of $g$ and  $\E_\infty g=0$, and because 
$E^\w_0(\si_{K(n)+1}-n)=E^\w_0E^\w_{X_n}(\si_1)$ is bounded.

%term is that the error above is bounded by the sum of three terms:
%two are of the kind
%\be
%\E\biggl[E^\w_0\biggl(\;\sum_{i=n}^{\si_{K(n)+1}-1} g(T_{X_i}\w)\biggr)
%E^\w_0\biggl(\;\sum_{j=0}^{n-1}  g(T_{X_j}\w)\biggr)\biggr].
%\label{error-1}
%\ee
%The first expression in the bracket is a bounded function of $\w$
%because $g$ is bounded and, by \eqref{4}, $E^\w_{X_n}(\si_1)$ is bounded uniformly
%in $\w$ and the state $X_n$.  By ergodicity of $\P_\infty$ the average
%$n^{-1} \sum_{j=0}^{n-1}  g(T_{X_j}\w)$ converges to $\E_\infty g=0$
%almost surely under the stationary environment chain. Hence
%this convergence happens almost surely  under $P^\w_0$ for $\P_\infty$-almost
%every  $\w$.  Then the same is true for  $\P$-almost
%every  $\w$. And finally by bounded convergence, the convergence to 0 is valid
%for the mean $E^\w_0\bigl(n^{-1} \sum_{j=0}^{n-1}  g(T_{X_j}\w)\bigr)$. Again this
%happens boundedly, so \eqref{error-1} is $\ord(n)$.
%
%The third error term is a product of two bounded terms such as
%the first expression in \eqref{error-1}. Hence, up to a negligible error,
%the variance of the quenched mean $E_0^\w(X_n)$ is equal to

Neglecting the error,  \eqref{error-1} equals 
\begin{align}
&\E\biggl[E^\w_0\biggl(\sum_x \sum_{i=0}^{\si_{K(n)+1}-1} \one\{ X_i=x\} g(T_x\w) \biggr)\nn\\
&\qquad\qquad\qquad 
\times \;
E^\w_0\biggl(\sum_y \sum_{j=0}^{\si_{K(n)+1}-1} \one\{ X_j=y\} g(T_y\w) \biggr)
\biggr].
\label{line-t3}
\end{align}
Manipulate each of the two expressions in the brackets as follows, 
noting that the
time spent at site $X_{\sigma_k}$ is $\si_{k+1}-\si_k$:
\begin{align*}
&E^\w_0\biggl(\sum_z \sum_{k=0}^{\si_{K(n)+1}-1} \one\{ X_k=z\} g(T_z\w)\biggr)\\
&\qquad=E^\w_0\biggl(\sum_z \sum_{k=0}^\infty \one\{ X_{\si_k}=z, \si_k\leq n\}
(\si_{k+1}-\si_k)  g(T_z\w)\biggr)\\
 &\qquad= \sum_z \sum_{k=0}^\infty P^\w_0(X_{\si_k}=z, \si_k\leq n)
E^\w_z(\si_{1})  g(T_z\w)\\
%&\qquad\qquad\qquad\qquad\text{[strong Markov property]} \\
&\qquad= \sum_z \sum_{k=0}^\infty P^\w_0(X_{\si_k}=z, \si_k\leq n)
  \frac{g(T_z\w)}{1-\pi_{zz}(\w)}.
\end{align*}
Combine the factors again, so that \eqref{error-1} turns into
\begin{align*}
&\E[|E_0^\w(X_n)-nv|^2]
= \ord(n) \\
&+\E\biggl[\sum_{x,y} \sum_{i,j\geq 0} P^\w_0(X_{\si_i}=x, \si_i\leq n)
 P^\w_0(X_{\si_j}=y, \si_j\leq n)
 \frac{g(T_x\w)}{1-\pi_{xx}(\w)}\, \frac{g(T_y\w)}{1-\pi_{yy}(\w)}\biggr]\,.
 \end{align*}
If, for example, $x<y$ then the factor $g(T_y\w)(1-\pi_{yy}(\w))^{-1}$ is
independent of the rest and is mean zero. Hence we simplify to 
\begin{align}
&\E[|E_0^\w(X_n)-nv|^2]
+ \ord(n) \nn\\
&=\E\biggl[\biggl|\frac{g(\w)}{1-\pi_{00}(\w)}\biggr|^2\biggr] \,
\E\biggl[\sum_x \sum_{i,j\geq 0} P^\w_0(X_{\si_i}=x, \si_i\leq n)
 P^\w_0(X_{\si_j}=x, \si_j\leq n)\biggr]\nn\\
&=\E\biggl[\biggl|\frac{g(\w)}{1-\pi_{00}(\w)}\biggr|^2\biggr] \,
E_{0,0}(\lvert X_{[0,n]}\cap \bar{X}_{[0,n]}\rvert),
\label{error-4}
\end{align}
where $X_I=\{X_i,i\in I\}$ and
$P_{0,0}=\E P_{0,0}^\w$ is the annealed process of two 
independent walks $X$ and $\bar X$
 in a common environment $\w$,  both starting at $0$.

Define now $L_0=0$ and for $j\geq1$,
\[L_j=\inf\{\ell>L_{j-1}: \exists n,m: X_n=\bar X_m=\ell\}.\]
Using similar arguments to those of Proposition \ref{hittingmoments} one
shows that $(L_j-L_{j-1})_{j\geq1}$ is an i.i.d.\ sequence under
$P_{0,0}$. Since these are non-negative
random variables, we have that $L_j/j$ converges to
$E_{0,0}(L_1)$.

Define $J_n=\max\{j:L_j\leq X_n\land \bar X_n\}$. Then
\[L_{J_n}\leq X_n\land\bar X_n\leq L_{J_n+1},\]
and Theorem \ref{lln}
implies that $J_n/n$ converges to $v/E_{0,0}(L_1)$. Since $J_n/n\in[0,1]$,
we have the same limit for $E_{0,0}(J_n)/n$. 
Dividing by $n$ in \eqref{error-4} and letting $n\to\infty$
 shows that the diffusion
coefficient is given by 
\[\kappa_m^2=\frac{v}{E_{0,0}(L_1)}\E[|E_0^\w(X_{\si_1}-v\si_1)|^2].\]
%For the convergence of the second moment of $n^{-1/2}E_0^\w(X_n)$ use the same
%argument as for the proof of part (c) of Proposition \ref{km+kq}.
\end{proof}

\subsection{The restricted-path case}
Now we consider the case where $d\geq 1$ is arbitrary but 
\be
\P(\,\exists z\not=0:\pi_{00}+\pi_{0z}=1)=1.
\label{the-case}
\ee
Before the theorem we go through some preliminaries.
Once an environment $\w$ has been fixed,
 each $x\in\Z^d$ has a unique point $w(x)\in\Z^d\setminus\{x\}$ such that
$\pi_{x,x}(\w)+\pi_{x,w(x)}(\w)=1$.
Once we assume forbidden direction and non-nestling (N)
with $\uhat$ in addition to \eqref{the-case},  any
$p>1$ in (M)  implies uniform boundedness of
the steps $w(x)-x$ and the existence of $\delta_1>0$
such that both $\pi_{x,w(x)}(\w)\geq \delta_1$
and $(w(x)-x)\cdot\uhat\geq\delta_1$
uniformly over $x$ and  $\P$-almost every $\w$.
It is convenient to replace the stopping times
$\{\si_j\}$ with the stopping times $\lambda_0=0$,
\[
\lambda_i=\inf\{n>\lambda_{i-1}: X_n\cdot\uhat\geq X_{\lambda_{i-1}}
\cdot\uhat +\delta_1\}\,,\qquad  i\geq 1.
\]
Replacing the increment  $1$ by $\delta_1$ in \eqref{sik}
does not affect the role of the stopping times in the proofs of 
the LLN or the annealed CLT. 
Hence we can express the asymptotic velocity as
$v=E_0(X_{\lambda_1})/E_0(\lambda_1)$ and the diffusion
matrix of the annealed invariance principle as
\be
\kD=\frac{E_0\bigl[(X_{\lambda_1}-\lambda_1v)
(X_{\lambda_1}-\lambda_1v)^t\bigr]}{E_0[\lambda_1]}
\label{ann-cov-la-1}
\ee
As in Theorem \ref{qclt-1d}, we will show next that the covariance
matrix \eqref{ann-cov-la-1} of the  annealed walk
can be decomposed as the limiting covariance of the
quenched walk plus the  limiting covariance of the
quenched mean process.
Noticing that in the case at hand
 $X_{\si_1}$ is not random under $P_0^\w$, define for the quenched mean
process the matrix $\kD_m$, analogously to \eqref{qmean-cov-la-1}, by
\be
\kD_m=\frac{\E\bigl[(X_{\lambda_1}-E^\w_0(\lambda_1)v)
(X_{\lambda_1}-E^\w_0(\lambda_1)v)^t\bigr]}{E_0[\lambda_1]}.
\label{qmean-cov-la-1-ii}
\ee
For the quenched walk, define the constant $\kappa_0$, analogously to \eqref{ell-var}, by
\[
\kappa_0^2=\frac{\E[\Var^\w_0(\lambda_1)]}{E_0[\lambda_1]}
=\biggl( \E\Bigl[\frac1{1-\pi_{00}}\Bigr]\biggr)^{-1}
\E\Bigl[\frac{\pi_{00}}{(1-\pi_{00})^2}\Bigr].
%\label{ell-var}
\]
The limit process for the quenched invariance principle 
 is $\kappa_0B(\cdot)v$ where $B(\cdot)$ is
a standard one-dimensional Brownian motion.
It has diffusion matrix
$\kD_q=\kappa_0^2vv^t$,  which is exactly what needs to be added to
\eqref{qmean-cov-la-1-ii} to yield \eqref{ann-cov-la-1}.

\begin{remark}
As it was  in the one-dimensional case, $\kappa_m>0$ if, and only if,
\begin{align}
\label{non-deg1-rep}
\P(D=v)=\P\Bigl(\w:w(0)= \frac{v}{1-\pi_{00}(\w)} \Bigr) <1.
\end{align}
If this condition fails, $Q_n^\w$ and $\Qtil^\w_n$ coincide.
Secondly,  $\kappa_q=0$ if, and only if, $P_0(\pi_{00}=0)=1$, i.e.\ \eqref{non-deg2} fails
to hold.
\end{remark}

\begin{theorem}
\label{hypE}
 Let $d\geq1$ and consider a product probability measure $\P$ on environments
with a forbidden direction $-\uhat\in\R^d\setminus\{0\}$
as in \eqref{forbidden}.
Assume non-nestling {\rm(N)} in direction $\uhat$,
the moment hypothesis {\rm(M)} with $p>1$, and {\rm\eqref{the-case}}.
\begin{enumerate}
\item For $\P$-almost every  $\w$ the
distributions $\Qtil_n^\w$ converge weakly to the distribution of
the process $\{\kappa_0 B(t)v : t\geq 0\}$
where $B(t)$ denotes standard one-dimensional Brownian motion.
\item Under $\P$, the scaled quenched mean process
$n^{-1/2}\{ E^\w_0(X_{[nt]})-ntv\}$ converges weakly to a Brownian
motion with diffusion matrix $\kD_m$.
\item Assume non-degeneracy condition \eqref{non-deg1-rep} holds.
%%if there are nonzero $x\neq z$ such that
%%both $\E\pi_{0x}>0$ and $\E\pi_{0z}>0$,
Then  for $\P$-almost every  $\w$, the sequence
of distributions $\{Q^\w_n:n\geq 1\}$ is not tight.
\end{enumerate}
\end{theorem}

\begin{proof}
The environment $\w$ uniquely determines a path
$\{0=z_0, z_1, z_2,  \dotsc\}$
such that $z_{i+1}=w(z_i)$. Each new point $z_i$  takes the walk to
 fresh environments. Consequently under $P_0$ the sequence
$\{\w_{z_i}\}$ of environments is also i.i.d. In particular,
$\{\xi_i=z_i-z_{i-1}\}$ are i.i.d.~random vectors with common distribution
$P_0\{\xi_i=z\}= \E[\pi_{0z}(1-\pi_{00})^{-1}]$ for $z\neq 0$
and mean $\tilde{v}= E_0\xi_i$.

Under $P^\w_0$ the walk $X_n$ is confined to the path
$\{z_i\}$.
Let $Y_n$ mark the location of the walk
along this path:
\[X_n=z_{Y_n}.\]
Given $\w$, $\{Y_n\}$ is a nearest-neighbor
 random walk on nonnegative integers
 with transitions
\[
P^\w_0(Y_{n+1}=i\vert Y_n=i)=\pi_{00}(T_{z_i}\w)=
1-P^\w_0(Y_{n+1}=i+1\vert Y_n=i).
\]
The previously defined $\lambda_i$ also give the hitting times
$\lambda_i=\inf\{n: Y_n=i\}$. (This is why we use them here.)  
Furthermore, under $P_0$ $X_{\lambda_1}$ has the same distribution
as $\xi_1=z_1$, and so
 $\tilde{v}= E_0\xi_i=E_0(\lambda_1)v$.

The limiting velocity for $Y_n$ is ($P_0$-a.s.)
$n^{-1}Y_n\to \gamma=(E_0\lambda_1)^{-1}$.
By Theorem \ref{qclt-1d} and Example \ref{1d-nn-example},
$Y_n$ satisfies a quenched invariance
principle
\be
\Bigl\{\frac{Y_{[nt]}-E^\w_0(Y_{[nt]})}{\sqrt{n}} : t\geq 0\Bigr\}
\mathop{\longrightarrow}^{\mathrm{dist.}}
\kappa_q B(\cdot)
\label{y-clt-2}
\ee
with limiting variance
\[\kappa_q^2=\gamma^3\E[\Var^\w_0(\lambda_1)]
=(E_0[\lambda_1])^{-2}\kappa_0^2. \]
Since the process $\kappa_q B(t)\vtil$ is the same as the limit
process
$\kappa_0 B(t)v$ claimed in part (a) of  the theorem,
part (a) will follow from  showing  that
\be
\lim_{n\to\infty} \;\max_{m\leq n}\; \frac1{\sqrt{n}} \Bigl\lvert
 \{X_{m}-E^\w_0(X_{m})\}
\;-\; \{Y_{m}-E^\w_0(Y_{m})\}\, \vtil\,\Bigr\rvert
=0
\quad\text{$P_0$-almost surely.}
\label{goal-1a}
\ee

Introduce the centered random vectors
  $\zbar_i=z_i-i\vtil$.
Decompose \eqref{goal-1a} into two tasks:
\be
\lim_{n\to\infty} \;\max_{m\leq n}\; \frac1{\sqrt{n}} \bigl\lvert
 \zbar_{Y_m}-\zbar_{[m\gamma]}
 \bigr\rvert
=0
\quad\text{$P_0$-a.s.}
\label{goal-1aa}
\ee
and
\be
\lim_{n\to\infty} \;\max_{m\leq n}\; \frac1{\sqrt{n}} \bigl\lvert
 \zbar_{[m\gamma]}-E^\w_0(X_{m})
+E^\w_0(Y_{m})\vtil\, \bigr\rvert
=0
\quad\text{$P_0$-a.s.}
\label{goal-1ab}
\ee

For \eqref{goal-1aa},
\begin{align}
&P_0 \Bigl( \max_{m\leq n}\;  \bigl\lvert
 \zbar_{Y_m}-\zbar_{[m\gamma]}
 \bigr\rvert \geq \e{\sqrt{n}} \Bigr)\nn\\
&\quad\leq \sum_{m=1}^n \Bigl\{  P_0(\,\lvert \,Y_m-[m\gamma]\,\rvert
\geq n^{3/4})\label{prob-11}\\
&\qquad\qquad
+ P_0 \bigl( \;\max_{k: \lvert k-[m\gamma]\rvert \leq n^{3/4}}\;  \bigl\lvert
 \zbar_{k}-\zbar_{[m\gamma]}  \bigr\rvert \geq \e{\sqrt{n}} \,\bigr)\Bigr\}.
\label{prob-12}
\end{align}

Both probabilities are controlled by standard large deviation
estimates.
For this we use  the following general inequality.

\begin{lemma}
Let $S_m=Z_1+\dotsm+Z_m$ be a sum of i.i.d.~mean zero
random variables  $\{Z_i\}$  with an exponential  moment:
$E(e^{\theta \abs{Z_1}})<\infty$ for some $\theta>0$.
Then there is a constant  $C$ determined by the
distribution of $\{Z_i\}$ such that, for all $m\geq 1$
 and $h>0$,
\be
P( \max_{k\leq m} S_k \geq h)\leq
e^{-C( \frac{h^2}{m}\land h)}. 
\label{ld-ineq-iid}
\ee
Same bound is of course valid for 
$\displaystyle{P( \min_{k\leq m} S_k \leq -h) =P( \max_{k\leq m} (-S_k) \geq h).}$
\label{ld-lm-1}
\end{lemma}

\begin{proof}[Proof of Lemma \ref{ld-lm-1}]
Let $\phi(\theta)=\log Ee^{\theta Z_1}$. Note that
$\phi(\theta)\geq 0$, and there exist constants $0<A,\theta_0<\infty$
such that
 $\phi(\theta) \leq A\theta^2$ for $\abs{\theta}\leq \theta_0$.
For $0\leq \theta\leq\theta_0$, applying Doob's inequality to the martingale
$e^{\theta S_k-k\phi(\theta)} $
gives
\[
P( \max_{k\leq m} S_k \geq h) \leq e^{ -\theta h +m \phi(\theta)}
\leq e^{ -\theta h +Am\theta^2}.
\]
If $h/(2Am)<\theta_0$, pick $\theta=h/(2Am)$ to get the bound
$e^{-h^2/(4Am)}$.  Now if $h/(2Am)\geq \theta_0$, then
\[
-\theta h +Am\theta^2\leq -\theta h +\theta^2h/(2\theta_0)
= -\theta_0h/2
\]
where we chose $\theta=\theta_0$. 
\end{proof}

The walk $Y_n$ is not a sum of i.i.d.~steps, but we can
consider the i.i.d.~sequence $\{\lambda_j-\lambda_{j-1}\}$
of successive sojourn times.  With $m_0=[m\gamma]+[n^{3/4}]$,
 $m_1=([m\gamma]-[n^{3/4}])^+$, and $b\in(0,1)$ a small positive
number, the probability on line \eqref{prob-11} is bounded as
follows:
\begin{align}
&P_0(\lvert Y_m-[m\gamma]\rvert \geq n^{3/4})\nn\\
&\qquad\leq
 P_0(\lambda_{m_0}\leq m)
+ P_0(\lambda_{m_1}\geq m)\one_{\{m_1>0\}}\nn\\
&\qquad\leq
 P_0(\lambda_{m_0}-E_0(\lambda_{m_0})\leq -bn^{3/4})
+ P_0(\lambda_{m_1} -E_0(\lambda_{m_1}) \geq bn^{3/4} )
\one_{\{m_1>0\}}\nn\\
&\qquad\leq e^{ -C(\frac{b n^{3/2}}{m_0}\land n^{3/4})}
+ e^{ -C(\frac{b n^{3/2}}{m_1}\land n^{3/4})}
\one_{\{m_1>0\}}\nn\\
&\qquad\leq 2e^{-Cn^{1/2}}.  \label{ld-bd-34}
\end{align}
A direct application of
 \eqref{ld-ineq-iid} to line \eqref{prob-12} gives a bound
$Ce^{-Cn^{1/4}}$. To summarize, we have the bound
\[
P_0 \Bigl( \max_{m\leq n}\;  \bigl\lvert
 \zbar_{Y_m}-\zbar_{[m\gamma]}
 \bigr\rvert \geq \e{\sqrt{n}} \Bigr)
\leq Cn\bigl( e^{-Cn^{1/2}}+ e^{-Cn^{1/4}}\bigr)
\]
which is summable in $n$. By Borel-Cantelli, \eqref{goal-1aa} has
been verified.

For \eqref{goal-1ab},  note that by the
 boundedness of $\xi_i-\vtil$,
$\lvert \zbar_{[m\gamma]}-\zbar_i\rvert\leq Cn$.
Write
\begin{align}
&\bigl\lvert \zbar_{[m\gamma]}-
E^\w_0(X_m)+E^\w_0(Y_m)\vtil \bigr\rvert
=\biggl\lvert
\sum_{i=0}^m P^\w_0(Y_m=i)(\zbar_{[m\gamma]}-\zbar_i)\biggr\rvert \nn\\
&\qquad
\leq  \sum_{i: \lvert i-[m\gamma]\,\rvert \,\geq\, n^{3/4}}
P^\w_0(Y_m=i)\bigl\lvert \zbar_{[m\gamma]}-\zbar_i\bigr\rvert \nn\\
&\qquad\qquad\qquad\qquad +
 \sum_{i: \lvert i-[m\gamma]\,\rvert \,<\, n^{3/4}}
P^\w_0(Y_m=i)\bigl\lvert \zbar_{[m\gamma]}-\zbar_i\bigr\rvert\nn\\
&\qquad
\leq Cn\sum_{i: \lvert i-[m\gamma]\,\rvert \,\geq\, n^{3/4}}
P^\w_0(Y_m=i)
\;+\;
 \max_{i: \lvert i-[m\gamma]\,\rvert \,<\, n^{3/4}}
\lvert \zbar_{[m\gamma]}-\zbar_i\rvert \nn\\
&\qquad
\leq Cn  P^\w_0\bigl(\,
\lvert Y_m-[m\gamma]\,\rvert \,\geq\, n^{3/4}\,  \bigr)
\;+\;
 \max_{i: \lvert i-[m\gamma]\,\rvert \,<\, n^{3/4}}
\lvert \zbar_{[m\gamma]}-\zbar_i\rvert. \nn
\end{align}
The last term is treated as was done above for \eqref{prob-12}.
For the second last term a summable deviation bound develops as follows.
\begin{align*}
&\P\Bigl( \max_{m\leq n} nP^\w_0\bigl(\,
\lvert Y_m-[m\gamma]\,\rvert \,\geq\, n^{3/4}\,  \bigr)
\geq \e\sqrt{n}\Bigr) \\
&\qquad\leq
\sum_{m=1}^n \frac{n}{\e\sqrt{n}}\,
P_0\bigl(\,
\lvert Y_m-[m\gamma]\,\rvert \,\geq\, n^{3/4}\,  \bigr)
\\
&\qquad\leq C n^{3/2} e^{-Cn^{1/2}}.
\end{align*}
Above we
repeated the  estimate from line \eqref{ld-bd-34}.
To summarize, we have the bound
\[
P_0 \Bigl( \max_{m\leq n}\;
\bigl\lvert \zbar_{[m\gamma]}-
E^\w_0(X_m)+E^\w_0(Y_m)\vtil \bigr\rvert
 \geq \e{\sqrt{n}} \Bigr)
\leq Cn^{3/2}\bigl( e^{-Cn^{1/2}}+ e^{-Cn^{1/4}}\bigr),
\]
and \eqref{goal-1ab} also has
been verified. This completes the proof of the convergence
of $\Qtil^\w_n$.

We turn to part (b) of the theorem.
To see the convergence of the quenched mean process,
note that
$nv=[n\gamma]\vtil+O(1)$ and  write
\be
\begin{split}
\frac{E^\w_0(X_n)-[n\gamma]\vtil}{\sqrt{n}}
&=\frac{E^\w_0(X_n)-E^\w_0(Y_n)\vtil-\zbar_{[n\gamma]}}{\sqrt{n}}\\
&\qquad 
 +\frac{\bigl(E^\w_0(Y_n)-[n\gamma]\bigr)\vtil+\zbar_{[n\gamma]}}{\sqrt{n}}.
\end{split}
\label{aux-res}
\ee
The first ratio after the equality sign
 tends to zero by \eqref{goal-1ab}.

Notice now that in \eqref{Delta}, applied to the $Y_n$ walk satisfying
$P_0(Y_1\in\{0,1\})=1$, we have $a_i=0$ for $i>0$, and
\[\Delta(\w)=\frac{E_0^\w(Y_1)-\gamma}{1-\pi_{00}(\w)}
=1-\gamma E_0^\w(\la_1).\]
 Therefore,
\eqref{ERn-small} and \eqref{EwXn} imply
%\be
%\begin{split}
%%\lim_{n\to\infty}\max_{k\leq n}\;\frac1{\sqrt{n}}\;
%%\lvert E^\w_0(X_k)-G_k(\w)\rvert \\
%%&
%\lim_{n\to\infty}\max_{k\leq n} \biggl\lvert
%\;\frac{E^\w_0(Y_k)-kv}{\sqrt{n}}\;-\;
%\frac1{\sqrt{n}}\sum_{j=0}^{[kv]-1} \bigl(1-\gamma E^{T_j\w}_0(\la_1)\bigr)
%\;\biggr\rvert=0
%\end{split}
%\label{1d-nn-qmean-1}
%\ee
%which shows
that on the diffusive scale the centered quenched mean process $E_0^\w(Y_n)-[n\gamma]$
behaves
like the sum of  i.i.d.\ mean zero random variables $1-\gamma E^{T_{z_j}\w}_0(\la_1)$,
$0\leq j\leq[n\gamma]-1$.
Expanding $\zbar_{[n\gamma]}$, we see that the
last term in \eqref{aux-res}  is asymptotically the same as
\be
\frac1{\sqrt{n}}\sum_{j=0}^{[n\gamma]-1}
\bigl(\xi_{j+1}-\gamma E^{T_{z_j}\w}_0(\lambda_1)\vtil\bigr).
\label{sum-35}
\ee
Under $\P$ this is a sum of mean zero i.i.d.~random vectors, and
this gives the invariance principle for $n^{-1/2}(E^\w_0(X_{[nt]})-ntv)$.

Lastly, part (c) is proved by comparison with (b).
Write
\begin{align*}
\frac{X_n-[n\gamma]\vtil}{\sqrt{n}}&=\frac{X_n-E^\w_0(X_n)}{\sqrt{n}}
+\frac{E^\w_0(X_n)-E^\w_0(Y_n)\vtil-\zbar_{[n\gamma]}}{\sqrt{n}}\\
&\qquad +\frac{\bigl(E^\w_0(Y_n)-[n\gamma]\bigr)\vtil+\zbar_{[n\gamma]}}{\sqrt{n}}.
\end{align*}
The first ratio after the equality sign is tight under $P^\w_0$
by part (a).  The second ratio tends to zero by \eqref{goal-1ab}.
As observed above, the
last ratio is asymptotically the same as
the sum \eqref{sum-35}
of mean zero i.i.d.~random vectors under $\P$.
 Hence under a fixed $\w$ this sequence of normalized sums
 cannot be bounded unless the summands are degenerate.
This is exactly the condition in the theorem.
\end{proof}

\begin{appendix}

\section{A renewal process bound}
\label{app-ren}
We write $\N^*=\{1,2,3,\dotsc\}$ and
$\N=\{0,1,2,\dotsc\}$.
The setting for the next technical lemma
is the following.
 $\{Y_i: i\in\N^*\}$ is a sequence of
 i.i.d.~positive
integer-valued random variables,
and $\{\Ytil_j: j\in\N^*\}$ is an independent copy of this
sequence.
 $Y$ will also denote a  random
variable with the same distribution as $Y_1$.
The corresponding renewal processes are defined by
\[
\text{$S_0=\Stil_0=0$, $S_n=Y_1+\dotsm+Y_n$
and $\Stil_n=\Ytil_1+\dotsm+\Ytil_n$ for $n\geq 1$. }
\]
Let $h$ be the largest positive integer such that
the common distribution
of $Y_i$ and $\Ytil_j$  is supported on $h\N^*$.
For $i,j\in h\N$ define
\[
L_{i,j}=\inf\{ \ell\geq 1: \text{ there exist $m,n\geq 0$ such that
 $i+S_m=\ell=j+\Stil_n$} \}.
\]
The restriction $\ell\geq 1$ in the  definition has the
consequence that $L_{i,i}=i$ for $i>0$ but $L_{0,0}$ is
nontrivial.

\begin{lemma}  Let $1\leq p<\infty$ be a real number, and assume
$E(Y^{p+1})<\infty$.  
  Then there
exists a finite constant $C$ such that
for all $i,j\in h\N$, 
\begin{align}
\label{firstgoal}
E(L_{i,j}^p) \leq C(1+i^p + j^p).
\end{align}
\label{renewal-lm}
\end{lemma}

\begin{proof} 
First some reductions.
If $h>1$, then we can consider the random variables
$\{h^{-1}Y_i, h^{-1}\Ytil_j\}$. Thus we may   assume
that $h=1$.

For convenience, we assume $i\leq j$ in $L_{i,j}$.
The cases $L_{i,i}$ for $i>0$ were already observed to be
trivial. If $i<j$ then
 $L_{i,j}$ has the same distribution as $i+L_{0,j-i}$.
 If we can prove the statement for
all $L_{0,j}$, then 
\[
E(L_{i,j}^p) \leq 2^p i^p + 2^p E(L_{0,j-i}^p)\leq
2^p i^p+2^p C(1+(j-i)^p) \leq 2^p(C+1)(1+i^p+j^p).
 \]
So we only need to consider the case $L_{0,j}$ for
$0\leq j<\infty$.

Next we remove the case $L_{0,0}$. Suppose we have
the  conclusion for $L_{0,j}$ for all $0<j<\infty$. 
On the event $\{Y_1=j\}$, $L_{0,0}=\tilde{L}_{j,0}$ computed
with the partially shifted
 sequences $(Y_2,Y_3,Y_4,\dotsc; \Ytil_1,\Ytil_2,\Ytil_3,\dotsc)$,
and furthermore $\tilde{L}_{j,0}$ is independent of $Y_1$
and  distributed like $L_{0,j}$.
Thus 
\be
\begin{split}
E(L_{0,0}^p)&=\sum_{j=1}^\infty E(Y_1=j, L_{0,0}^p) 
=\sum_{j=1}^\infty E(Y_1=j, \tilde{L}_{j,0}^p)
= \sum_{j=1}^\infty P(Y_1=j) E(L_{0,j}^p)\\
&\leq 
C\sum_{j=1}^\infty P(Y_1=j) (1+j^p) = C(1+ E[Y^p])<\infty. 
\end{split}\label{L00-arg}
\ee

Now we turn to proving the statement \eqref{firstgoal} 
for $L_{0,j}$, $0<j<\infty$.  
The proof will make use of several Markov chains. 

Denote the forward recurrence time  of the
pure renewal process $\{S_n:n\geq 0\}$ by
\[
B^0_k=\min\{ m\geq 0:  k+m\in\{S_n:n\geq 0\}\}.
\]
More  generally, if we start the renewal process with a deterministic
delay $Y_0=r\in\N$ and then define
$S_n=Y_0+Y_1+\dotsm+Y_n$ for $n\geq 0$,
 we denote the forward recurrence time
by $B^r_k$, defined exactly as above. $\{B^r_k: k\geq 0\}$ is  a
Markov chain with initial state $B^r_0=r$ and transition
probability
\[
\text{$p(0,j)=P(Y=j+1)$ for $j\geq 0$, and $p(j,j-1)=1$ for $j\geq 1$.}
\]
The usefulness of this Markov chain is that $B^r_k=0$ if, and only if,
$k\in\{S_n:n\geq 0\}$.

The state space for the  $p$-chain is
$\Gamma=[0,\rho)\cap\N$ where 
\[
\rho=\sup\{ y : P[Y=y]>0\} \in\N^*\cup\{\infty\}
\]
is the supremum of the support of $Y$.
This chain is irreducible, and the assumption $h=1$
makes it aperiodic. Its unique
 invariant distribution is given by
\[
\pi(j)=\frac{P(Y>j)}{E(Y)}.
\]
We also consider the joint chain $(B^r_k,\Btil^s_k)$
of two independent  copies
of the Markov chain $B^r_k$, started at $(r,s)\in\Gamma^2$,
 with transition
\[
\pbar((x,y),(u,v))=p(x,u)p(y,v).
\]
The familiar argument from Markov chain theory shows that
the irreducibility and aperiodicity of transition $p$
implies the irreducibility of  transition $\pbar$;
see \cite[p. 129-130]{resnick} or \cite[p. 311]{durrett}.
The unique invariant distribution of $\pbar$ is $\pibar(i,j)=\pi(i)\pi(j)$.
We are ready to prove \eqref{firstgoal} for bounded $Y$.

\medskip

{\bf Case 1.} $\rho<\infty$.

\medskip

For $0\leq i,j<\rho$,
$L_{i,j}$ has the same distribution as
the first return time
\[
T_{(0,0)} =\inf\{ k\geq 1: (B^i_k,\Btil^j_k)=(0,0)\}
\]
to state $(0,0)$, given that the joint chain starts from
$(i,j)$. The reason is that considering $i+S_m$ and $j+\Stil_n$
 is the same
as starting the renewal processes with  delays $i$ and $j$.
Since the state space of the $\pbar$-chain is a single positive
recurrent class,  we can pick for each state  $(r,s)\in\Gamma^2$
a finite path to $(0,0)$ of positive probability. Since there are
only finitely many states, 
there exist $n_0<\infty$ and $\delta>0$ such that 
\[
\inf_{(r,s)\in\Gamma^2} P_{(r,s)}[ T_{(0,0)}\leq n_0 ] \geq \delta.
\]
>From this and the Markov property follow uniform exponential tail
bounds
\[
\sup_{(r,s)\in\Gamma^2} P_{(r,s)}[ T_{(0,0)}\geq n ] \leq C\eta^n 
\]
for constants $0< C<\infty$ and $0<\eta<1$, and thereby 
\[
A_0\equiv \sup_{0\leq j<\rho} E(L_{0,j}^p) <\infty.
\]
To cover $j\geq \rho$  we can repeat the argument given above
in \eqref{L00-arg}.  Let $U$ be the stopping time for the sequence
$(Y_i)$  defined
by $S_{U-1}<j\leq S_U$. On the event $S_U-j=B^0_j=x>0$, let $\tilde{L}_{x,0}$
be defined  in terms of the partially restarted sequences
$(Y_{U+1}, Y_{U+2}, Y_{U+3},\dotsc; \Ytil_1,\Ytil_2,\Ytil_3,\dotsc)$.
$\tilde{L}_{x,0}$ is independent of $B^0_j$ and distributed like 
$L_{0,x}$.  We have $L_{0,j}= j+\tilde{L}_{x,0}$ because $\tilde{L}_{x,0}$
ignores the distance from $0$ to $j$.
\begin{align*}
E(L_{0,j}^p)&=j^p P(B^0_j=0)+\sum_{0<x<\rho} E\bigl[B^0_j=x, (j+\tilde{L}_{x,0})^p\,
\bigr] \\
&\leq 2^pj^p + 2^p \sum_{0<x<\rho}  P(B^0_j=x) E(L_{0,x}^p)
\leq  2^pj^p + A_0.
\end{align*}

\medskip

{\bf Case 2.} $\rho=\infty$.

\medskip

To handle this case 
we define another Markov chain on non-negative integers
with transition $q$.
The chain absorbs at zero so $q(0,0)=1$.
The other transitions are
\[
q(x,y)=P( B^0_x=y),\text{ for $x\geq 1$ and $y\geq 0$.}
\]
Let $\{\zeta_i:i\geq 0\}$ denote the process with transition $q$.
This chain represents a construction of  $L_{0,j}$ where 
we take turns sampling variables $Y_i$ and $\Ytil_k$ 
until the sums $S_m$ and $\Stil_n$  meet. 
  Assume now that $0<j<\infty$. For $j=0$
the description below would stop right away and not work correctly.

(i) Put $\zeta_0=j$.

(ii)  First, construct $0=S_0< S_1<S_2<\dotsm$
until the first $S_m\geq j$. If this $S_m=j$ we are done and
we set $\zeta_1=0$ and $L_{0,j}=S_m=j$. If $S_m$ overshoots so that $S_m>j$, then let
$\zeta_1=S_m-j=B^0_j$ be the forward recurrence time
at $j=\zeta_0$.

(iii) Now turn to construct $0=\Stil_0<\Stil_1<\Stil_2<\dotsm$
until the first  $\Stil_n\geq \zeta_1$. If $\Stil_n=\zeta_1$
 we are done, and then $\zeta_2=0$ and  $L_{0,j}=\zeta_0+\zeta_1$.
If $\Stil_n>\zeta_1$ overshoots, then let
$\zeta_2=\Stil_n-\zeta_1=\Btil^0_{\zeta_1}$  be
the new forward recurrence time.

(iv) Next go back to add more terms to the first
sum $S_m$. The earlier terms do not play a role. Given
$\zeta_2$, we add  new independent terms  $Y'_1, Y'_2, Y'_3, \dotsc$
until  $Y'_1+\dotsc+Y'_m\geq\zeta_2$.
And the overshoot is then $\zeta_3=Y'_1+\dotsc+Y'_m-\zeta_2=
B^0_{\zeta_2}$,
a forward recurrence time  relative to a new renewal
process independent of the earlier ones.

These  steps are repeated until
 the sums  meet. This happens at time
\[
\nu_0=\inf\{i\geq 1: \zeta_i=0\}
\]
and then 
\be
\text{for $0<j<\infty$,} \quad 
L_{0,j}= \sum_{i=0}^{\nu_0}\zeta_i \quad\text{with $\zeta_0=j$.}
\label{L0j-aux}
\ee
 The
term $\zeta_{\nu_0}$ is zero but for conditioning it is
 convenient  to extend the sum to $\nu_0$.

To make use of this we need two auxiliary bounds: a uniform
bound on $E_x(\zeta_1^p)$, and a tail bound on $\nu_0$.

\begin{lemma} Let $1\leq p<\infty$, a real number. Then
for some constant $C$ independent of $p$, 
$E_x(\zeta_1^p)\leq C E(Y^{p+1})$  for all $x\geq 0$.
\label{renewal-aux-lm-1}
\end{lemma}

\begin{lemma} Assume $\rho=\infty$. 
 Then there exist $0<\theta<1$ and $A<\infty$
such that $P_x(\nu_0>n)\leq A\theta^n$   for all $x\geq 1$
and $n\geq 0$.
\label{renewal-aux-lm-2}
\end{lemma}

Let us complete the proof of the main statement and then
return to prove the auxiliary lemmas. Recall that  $0<j<\infty$.
 As all the summands
are non-negative, we can manipulate the sums without regard
to finiteness.
\[
\begin{split}
E(L_{0,j}^p)&= E_j\biggl[\,\biggl(\sum_{i=0}^{\nu_0}\zeta_i\biggr)^p\;\biggr] 
= E_j\biggl[\,\biggl(j+\sum_{i=1}^{\nu_0}\zeta_i\biggr)^p\;\biggr] \\
&\leq 2^pj^p \;+\; 2^p E_j\biggl[\,\biggl(\sum_{i=1}^{\nu_0}\zeta_i\biggr)^p\;\biggr].
\end{split}
%\label{L0j-aux-2}
\]
It remains to show that the last expectation is bounded by a constant
independently of $j$.  
Apply Minkowski's inequality (with a limit to infinitely many terms), 
note that the event $\{\nu_0\geq i\}$
is measurable with respect to
$\sigma\{\zeta_0,\dotsc,\zeta_{i-1}\}$, and apply Lemmas 
  \ref{renewal-aux-lm-1} and  \ref{renewal-aux-lm-2}.
%Let $I$ denote the set of vectors
%$(i_1,\dotsc,i_p)\in (\N^*)^p$ whose coordinates are not all equal. 
\begin{align*}
&\biggl\{ 
E_j\biggl[\,\biggl(\sum_{i=1}^{\infty}\ind_{\{\nu_0\geq i\}}\zeta_i\biggr)^p\;\biggr]
\,\biggr\}^{1/p} 
\leq \sum_{i=1}^{\infty} \bigl\{ \,E_j(\ind_{\{\nu_0\geq i\}}\zeta_i^p) \,\bigr\}^{1/p}\\
&= \sum_{i=1}^{\infty} \bigl\{ \, E_j\bigl(\ind_{\{\nu_0\geq i\}}E_{\zeta_{i-1}}[\zeta_1^p]\bigr) 
\,\bigr\}^{1/p}\\
&\leq \sum_{i=1}^{\infty} C\{E(Y^{p+1})\}^{1/p} \bigl\{ P_j(\nu_0\geq i) \bigr\}^{1/p}
 \leq C_1 < \infty,
\end{align*} 
 bounded independently of $j$ as required. This completes the proof of Lemma
\ref{renewal-lm}. 
\end{proof}

It remains to prove the auxiliary lemmas used above.

\begin{proof}[Proof of Lemma \ref{renewal-aux-lm-1}]
This lemma is trivial for $x=0$ because the $\zeta$-chain
absorbs at $0$.  For $x\geq 1$ a 
straightforward calculation, utilizing $S_n\geq n$, gives
\begin{align*}
E_x(\zeta_1^p)&\leq  C\,\sum_{k=1}^\infty P(B^0_x\geq k)k^{p-1}
=C\,\sum_{k=1}^\infty \sum_{n=1}^x \sum_{i=1}^x
P(S_{n-1}=x-i, S_n\geq x+k )k^{p-1}\\
&=C\,\sum_{k=1}^\infty \sum_{n=1}^x \sum_{i=1}^x
P(S_{n-1}=x-i) P( Y\geq i+k )k^{p-1}\\
&=C\,\sum_{k=1}^\infty  \sum_{i=1}^x
P(\text{$S_{n}=x-i$ for some $n$}) \sum_{m=i+k}^\infty
 P( Y= m ) k^{p-1} \\
&\leq C\,\sum_{k=1}^\infty  \sum_{i=1}^x \sum_{m=i+k}^\infty
 P( Y= m )k^{p-1}\\
&\leq  C\,\sum_{m=2}^\infty P( Y= m ) \sum_{k=1}^\infty  \sum_{i=1}^\infty
\one\{i+k\leq m\} k^{p-1}\\
&\leq C\,E[Y^{p+1}].
\end{align*}
The third last inequality came from bounding probabilities by $1$,
the  second last inequality from removing the upper bound $x$
from the index $i$.
\end{proof}

\begin{proof}[Proof of Lemma \ref{renewal-aux-lm-2}]
Since the transition $p$ of $B^r_k$ is
positive recurrent, irreducible  and aperiodic,
\[
\lim_{x\to\infty} P(B^0_x=0)
=\lim_{k\to\infty} p^k(0,0)=\pi(0)>0.
\]
Fix $0<\e_0<1$ and $m_0<\infty$ so that
$P(B^0_x=0)\geq \e_0$ for $x\geq m_0$.

Now we use the assumption $\rho=\infty$.
Since the support of the distribution of $Y$ is unbounded,
\[
\e_1=P(Y\geq 2m_0)>0.
\]
Then for $0<x\leq m_0$
\[
P(B^0_x\geq m_0) \geq P( Y_1\geq x+m_0)\geq \e_1.
\]
In terms of the transition $q$,
for $x<m_0$,
\[
q^2(x,0)\geq \sum_{y\geq m_0} q(x,y)q(y,0) =
 \sum_{y\geq m_0} P(B^0_x=y)P(B^0_y=0) \geq \e_0\e_1,
\]
while for $x\geq m_0$,
\[
q(x,0)=P(B^0_x=0)\geq \e_0\geq \e_0\e_1.
\]

From this follows, for any $n\geq2$, $x\geq 1$,    and with
$\gamma=1-\e_0\e_1$,
\begin{align*}
P_x(\nu_0>n)&=\sum_{z\geq 1}P_x(\nu_0>n-2, \zeta_{n-2}=z)
P_z(\zeta_1\neq 0, \zeta_2\neq 0)\\
&\leq \sum_{1\leq z<m_0}P_x(\nu_0>n-2, \zeta_{n-2}=z)P_z(\zeta_2\neq 0)\\
&\qquad\qquad +
 \sum_{ z\geq m_0}P_x(\nu_0>n-2, \zeta_{n-2}=z)P_z(\zeta_1\neq 0)\\
&\leq \gamma P_x(\nu_0>n-2)
\leq \dotsm\leq \gamma^{[n/2]}.
\end{align*}
\end{proof}

\end{appendix}

\begin{acknowledgment}
A large part of this work was done while
Firas Rassoul-Agha was at the Mathematical Biosciences Institute, Ohio State
University.
\end{acknowledgment}

\bibliographystyle{abbrv} % A regular style that is supported 
%\bibliographystyle{aop} % An AOP-like style found and edited by Firas
                        % (the tabbing is not right though)
\bibliography{refsfiras}

%% If you want to use thebibliography look below:

%% Uncomment this if you want to use Bibliography instead of References
%\renewcommand{\refname}{Bibliography}

%% Uncomment this, if you wanted the references in small font.
%% Don't forget to uncomment the end part!
%\begin{footnotesize} %or \begin{small}

%\begin{thebibliography}{99}
%Separate the author, title, journal, issue, pages, and year.
%It will make it easier later to put them in different formats.
%Example:
%
%\bibitem{ferrarifontes}
%P.A.~Ferrari and L.R.G.~Fontes:
%Title.....
%Electron. J. Probab.
%\textbf{3},
%1--34
%(1998)
%

%\end{thebibliography}

%\end{footnotesize} %or \end{small}

\end{document}